\documentclass[10pt,a4paper,oneside,reqno]{amsart}

\usepackage{amsmath,latexsym}
\usepackage{amssymb}
\usepackage{geometry}
\usepackage{hyperref}
\usepackage{multicol}
\usepackage{multirow}
\usepackage{float}
\usepackage{caption, subcaption}
\usepackage{algorithmic,algorithm}
\usepackage{adjustbox}

\newcommand{\dsum}{\displaystyle\sum}
\newcommand{\dmin}{\displaystyle\min}
\newcommand{\dmax}{\displaystyle\max}
\def\R{\mathbb{R}}

\def\E{\mathbb{E}}
\def\A{\mathcal{A}}
\def\U{\mathcal{U}}
\def\P{\mathbb{P}}
\def\x{\mathbf{x}}
\newtheorem{thm}{Theorem}
\newtheorem{lem}{Lemma}

\newtheorem{prop}{Proposition}
\newtheorem{cor}{Corollary}

\usepackage{tikz}
\usepackage{pgfplots}
\usepgfplotslibrary{groupplots}
\usetikzlibrary{shapes.geometric, arrows,calc}
\usepackage{color}

\let\origmaketitle\maketitle
\def\maketitle{
  \begingroup
  \def\uppercasenonmath##1{} % this disables uppercasing title
  \let\MakeUppercase\relax % this disables uppercasing authors
  \origmaketitle
  \endgroup
}

\begin{document}

\title[Minimal Radius Enclosing Polyellipsoids]{\huge Minimal Radius Enclosing Polyellipsoids}

\author[V. Blanco \MakeLowercase{and} J. Puerto]{{\large V\'ictor Blanco$^\dagger$ and  Justo Puerto$^\ddagger$}\medskip\\
$^\dagger$IEMath-GR, Universidad de Granada\\
$^\ddagger$IMUS, Universidad de Sevilla}

\address{IEMath-GR, Universidad de Granada, SPAIN.}
\email{vblanco@ugr.es}

\address{IMUS, Universidad de Sevilla, SPAIN.}
\email{puerto@us.es}

\thanks{This research has been partially supported by Spanish Ministry of Econom{\'\i}a and  Competitividad/FEDER grants number
MTM2016-74983-C02-01.}

\date{\today}

\begin{abstract}
In this paper we analyze  the extension of the classical smallest enclosing disk problem to the case of the location of a polyellipsoid to fully cover a set of demand points in $\mathbb{R}^d$. We prove that the problem is polynomially solvable in fixed dimension and analyze mathematical programming formulations for it. We also consider some geometric approaches for the problem in case the foci of the polyellipsoids are known. Extensions of the classical algorithm by Elzinga-Hearn are also derived for this new problem. Moreover, we address several extensions of the problem, as the case where the foci of the enclosing polyellipsoid are not given and have to be determined among a potential set of points or the induced covering problems when instead of polyellipsoids, one uses ordered median polyellipsoids. For these problems we also present Mixed Integer (Non) Linear Programming strategies that lead to efficient ways to solve it.
\end{abstract}

\maketitle

\begin{quotation}
\centerline{\textbf{\em This article is dedicated to Prof. Marco A. L\'opez on the occasion of his 70th birthday.}}
\end{quotation}

\section{Introduction}

Given a sets of demand points in a given finite dimensional normed space $\R^d$, Continuous Facility Location Problems (CFLP) deal with the determination of the optimal placement of some new points in order to minimize certain measure of the distances from the given to the new points. The most popular CFLP is the Weber Problem \cite{weber} in which a single facility is to be located minimizing the overall (weighted) sum of the Euclidean distances from the demand points to the facility. A common approach to solve the Weber Problem is via the Weiszfeld algorithm \cite{weiszfeld}. The extension to multiple facilities may either consider that the overall sum from the demand points to all facilities (multiple-allocation case) or to its closest facility (single-allocation case) has to be minimized. Different \textit{aggregation} measures, apart from the overall sum, have been also proposed to find optimal facilities for CFLP (see e.g., \cite{BPE13}).

The level curves of the sum of Euclidean distances to given points on the plane (the Euclidean planar Weber level curve), in case the ($n$) demand points (also known as foci within this framework) are not collinear, are called in the literature \textit{polyellipses}, \textit{$n$-ellipses} or \textit{multifocal ellipses}.  The regions bounded by those curves are known to be convex bodies and can be also defined in higher dimensional spaces and with different norm measures, inducing the notion of \emph{polyellipsoids}. Several authors have been interested in the analysis of these convex bodies, from Maxwell \cite{jcm} to Er\"dos and Vincze \cite{erdos}, mainly from an algebraic or geometric point of view (see \cite{parrilo}). In particular,  the problem whether all the convex plane curves can be arbitrarily approximated by polyellipses under a sufficiently large number of points (also known as \textit{foci}) was posed by Weiszfeld. In this sense, the Erd\"os-Vincze's theorem states that regular triangles cannot be presented as the Hausdorff limit of polyellipses even if the number of the foci can be arbitrary large. Recently, in \cite{vincze2018}, the authors extend the analysis for regular polygons in the plane. However, these sets (polyellipsoids) have only been partially analyzed under optimization lens beyond their connections with the Weber problem (see, e.g. \cite{polyellipses}).

On the other hand, Covering Location Problems consist of locating facilities to satisfy the demand of a set of customers whenever they are only allowed to be served within certain coverage areas. In the continuous case there are two main families of problems in this field: Full Covering Location problems (FCP) and Maximal Covering Location problems (MCP). In the former all the demands points have to be covered and the goal is to minimize the set up costs, which may consists of the cost of opening a facility and/or the cost of enlarging the coverage area of the facilities to reach the demand points (see \cite{berman2010,plastria}). In the latter (MCP), the goal is to locate the facilities such that the (weighted) number of points minus the facilities set-up costs is maximized. Both problems have been analyzed in the literature, mainly on the plane. If the number of facilities to open $p$ is given and the regions are defined as distance-based balls, the FCP is also called the continuous $p$-center  problem or the minimum volume enclosing sphere, since the problem can be interpreted as the one of locating $p$ facilities in the space, such that the maximum distance from the demand points to its closest facility (center) is minimized. The particular case in which a single facility is to be located and the distance measure is the Euclidean, Elzinga and Hearn provided one of the most popular algorithms in continuous location for this problem~\cite{EH}. More recently, linear time algorithms have been proposed for the weighted minimum sphere enclosing problems (\cite{Dyer86} and \cite{Megiddo83}). Also, the problem of finding the minimum volume enclosing ellipsoid in fixed dimension is known to be solved in linear time \cite{Dyer}. On the other hand, for the FCP in case of locating a single circle-shaped facility on the plane there is always an optimal solution of the problem whose position lies either in the intersection of certain circles centered at the demand points, or at the demand points, the so-called \textit{circle intersection points} (CIP) \cite{CL86,Drezner81}. This property allows one to transform the original continuous problem into a discrete covering location problem. The multifacility case is addressed in \cite{Church84}, proving the same CIP property, for both the Euclidean and the $\ell_1$-norm distance measures. In \cite{AM13} and \cite{CM09} the case of locating elliptical regions in the plane is studied, by means of Mixed Integer Non Linear Programming formulations and heuristic approaches. Apart from that, as far as we know, there are no further advances on the topic. In particular, the definition and use of weighted  $d$-dimensional polyellipsoids to fully cover a set of demand points have not been previously considered.

In this paper, we extend the classical smallest enclosing disk problem (also known as \textit{minimax facility location problem}) to the case in which instead of spheres or ellipsoids, polyellipsoids (higher dimensional polyellipses) are to be located. Since polyellipses with a single focus coincide with disks, our approach naturally extends the classical problems.

\begin{figure}[h]
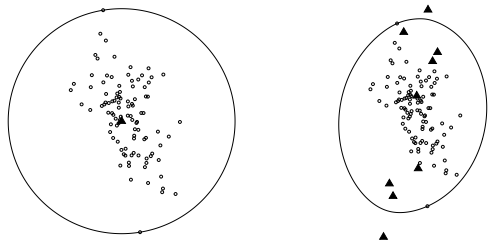

\begin{center}
\input{graph1.tex}~\input{graph2.tex}
\end{center}
\caption{Minimum radius circle (left) and $9$-ellipse (right) enclosing a set of demand points.\label{fig:intro}}
\end{figure}

Turning to the applications of these models, it is natural to think about locating several alarm sirens or facilities (which play the role of the foci) from where some emergency services should be dispatched. In these situations the service quality decreases with distance, so the closer the better. Therefore, the problem can be seen as to adequately locate a set of facilities whose relative positions are given and determine the minimum coverage distance for which the demand of all customers is satisfied. Note that the coverage area in this problem is defined as the region for which the sum of the distances from a customer to all the foci is less than or equal to the given threshold. In Figure \ref{fig:intro}-right we draw a situation in which the demand of a set of customers ($\circ$) has to be covered by nine facilities ($\blacktriangle$). The minimum radius coverage region is limited by the drawn polyellipse. As opposed to the case where the coverage has to be done with respect to a single focus (circle) which is shown in Figure \ref{fig:intro}-left.
Another field of applications, different from Location Analysis, comes from Data Analysis and outliers detection. Minimum enclosing balls allow one to determine possible outliers in a dataset, looking at the points that do not fall inside the ball. However, many clouds of points do not fit well to a sphere but to a more complex shape. Thus, polyellipsoids with an adequate set of foci allow one a more appropriate  fitting and a better outlier detection tool. A more sophisticated tool to detect outliers (and also for one-class classification) that is becoming popular in Data Science is Support Vector Data Description (SVDD) \cite{tax04}. In SVDD, outlier observations are detected by enclosing part of the points inside a Euclidean ball, after introducing misclassifying errors and transforming the points into a higher dimensional space using kernels. Our methodology could be adapted to this framework determining outliers observations for non-spherical data, avoiding the use of kernel transformations.

The paper is organized in seven sections. In Section \ref{sec:2} we introduce the framework of the paper, recall the definitions and some preliminary results of polyellipsoids and formulate the problem of determining the optimal position of a given set of foci such that a set of demands points is fully covered by the polyellipse with minimum radius. We also derive a polynomial-time complexity result for the problem. In Section \ref{sec:3} we describe several solution approaches for the problems, based either on dualization of the problem or in decomposition approaches inspired on the Elzinga-Hearn algorithm for the Euclidean $1$-center problem. Section \ref{sec:5} is devoted to report the results of an extensive computational experience which shows the efficiency of our decomposition approach. We analyze in Section \ref{sec:4} the very particular case in which the demand points are one-dimensional and derive a linear time algorithm (in the number of foci) for the minimal enclosing one-dimensional polyellipsoid. In Section \ref{sec:6} we show two possible extensions for the problem. First we study the problem of selecting the foci from a potential set of candidates when they are not \textit{``a priori''} given. Second, we present some results generalizing the notion of polyellipsoid, and the problem under analysis, to the ordered median framework. Finally, we draw some conclusions and point out some further research topics in Section \ref{sec:7}.

\section{Minimal enclosing polyellipsoids with given foci\label{sec:2}}

In this section we introduce the notation and concepts to be used for the rest of the paper. We recall the definition of polyellipsoid and derive some of its useful properties. We also analyze the problem of locating a polyellipsoid with given relative positions of the foci (but unknown radius and position) to minimally cover a given set of demand points.

\begin{figure}[h]
\begin{center}
\input{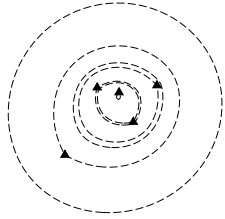}
\end{center}
\caption{Polyellipses with five foci and different radii.\label{fig00}}
\end{figure}

We are given a finite set of points $\U=\{u_1,\ldots,u_k\}$ in $\R^d$ and weights ${\bf \omega} \in \R^{k}_+$, and a distance measure induced by a norm $\|\cdot\|$. Without loss of generality, we assume that $\sum_{u\in \U} \omega_ u =1$. The \textit{minisum} or Weber location problem consists of finding the placement of a facility $x \in \R^d$ that minimizes the ${\bf \omega}$-weighted distances from $x$ to the points in $\U$, i.e., the function:
$$
\Phi_{\U,\mathbf{\omega}}(x) = \dsum_{u \in \U} \omega_u \|u-x\|, x \in \R^d.
$$
If the points in $\U$ are not collinear and the norm is strictly convex, $\Phi_{\U,\mathbf{\omega}}$ is strictly convex, and therefore has a unique minimum. The levels curves of $\Phi_{\U,\mathbf{\omega}}$, for $r\ge 0$, are given by the following regions:
$$
\mathbb{E}_{\U,\mathbf{\omega}}(r)= \Big\{x \in \R^d:  \dsum_{u \in \U} \omega_u \|u-x\| = r\Big\}.
$$
The set $\mathbb{E}_{\U,\mathbf{\omega}}(r)$ is called a (weighted) polyellipsoid with foci $\U$, weights ${\bf \omega}$, and radius $r$. Equal-weights Euclidean polyellipsoids on the plane are called \textit{polyellipses}. Polyellipses were introduced for the first time by Tschirnhaus in 1686. In case the number of foci is exactly one, polyellipses coincide with circles, and for two foci, one obtains standard ellipses. In Figure \ref{fig00} we draw some planar polyellipses with five foci ($\blacktriangle$) and different radii.

Let us denote by  $\mathbb{P}_{\U,\mathbf{\omega}}(r)$ the region bounded by the polyellipsoid $\mathbb{E}_{\U,{\bf \omega}}(r)$, i.e.,
$$
\mathbb{P}_{\U,\mathbf{\omega}}(r)= \Big\{x \in \R^d:  \dsum_{u \in \U} \omega_u \|x-u\| \leq r\Big\}.
$$
Abusing of notation, we will also call the sets $\mathbb{P}_{\U,\mathbf{\omega}}(r)$ polyellipsoids.

In the following we recall some known properties of the regions $\mathbb{P}_{\U,\mathbf{\omega}}(r)$:
\begin{prop}\label{prop:1}
Let $\U=\{u_1, \ldots, u_k\} \subseteq \R^d$, $\omega \in \R^{k}_+$ with $\sum_{u \in \U} \omega_u=1$ and $r\geq 0$. Then:
\begin{enumerate}
\item $\mathbb{P}_{\U,\mathbf{\omega}}(r)$ is a  closed, bounded and convex set.
\item If $\x^* \in \arg\dmin_{x \in \R^d} \sum_{u \in \U} \omega_u \| u- {\bf x} \|$ and $r^*=\sum_{u \in \U} \omega_u \| u- \x^* \|$:
\begin{enumerate}
\item $\x^* \in \mathbb{P}_{\U,\mathbf{\omega}}(r^*)$. If $\|\cdot \|$ is strictly convex and $\U$ are not collinear then $\mathbb{P}_{\U,\mathbf{\omega}}(r^*)= \{\x^*\}$.
\item $\mathbb{P}_{\U,\mathbf{\omega}}(r)= \varnothing$ for all $r<r^*$.
\item If $\mathbb{P}_{\U,\mathbf{\omega}}(r)=\{\bar \x\}$ for some $r>0$ and $\bar x\in \R^d$, then $r=r^*$ and $\bar \x=\x^*$.
\item $\x^* \in \mathbb{P}_{\U,\mathbf{\omega}}(r)$ for all $r\geq r^*$.
\end{enumerate}
\item If $d=2$ and $\|\cdot\|$ is the Euclidean norm. Then, $\mathbb{P}_{\U,\mathbf{\omega}}(r)$ can be represented as the set of solutions of a polynomial equation with degree $2k$ (if $k$ is odd) or $2k$ - $k \choose {\frac{k}{2}}$ (if $k$ is even) which can be expressed as the determinant of a symmetric matrix of linear polynomials~\cite{parrilo}.
\item $\mathbb{P}_{\U,{\bf \omega}}(r)$ is contained in the ring with center at $\x^*$ and radii $\dfrac{r-r^*}{\sum_{u\in \U} \omega_u}$ and $\dfrac{r+r^*}{\sum_{u\in \U}  \omega_u}$
\end{enumerate}
\end{prop}
The point $\x^*$ of the above result is known as the center of the polyellipsoid $\mathbb{P}_{\U,\mathbf{\omega}}(r^*)$ and it can be efficiently computed with classical location analysis algorithms as the hyperbolic modification of Weiszfeld algorithm \cite{weiszfeld}.

Although the above result gives us information about the complexity of determining the points which belong to a given polyellipsoid in terms of the number of foci, if the considered norm is a $\ell_p$-norm with $p\in \mathbb{Q}$, $p\ge 1$, these convex bodies are not only convex sets, but Second Order Cone (SOC) representable sets (see e.g.\cite{BPE14}).  The main implications of such a representation is that optimization problems with linear constraints and SOC-representable feasible sets can be efficiently solved via interior point methods~\cite{NN94} available in most common optimization solvers.

\subsection{The Minimum Radius Enclosing Pollyellipsoid Problem with Given Foci}
Given a finite set of demand points, $\A=\{a_1,\ldots,a_n\} \subseteq \R^d$, the goal of the \textit{Minimum Radius Enclosing Pollyellipsoid Problem with Given Foci} (MEP, for short) is to determine the placement and the minimum radius of the polyellipsoid with foci $\mathcal{U}=\{u_1,\ldots, u_k\}$ and weights $\mathbf{\omega}\in \R^k_+$ with $\sum_{j=1}^k \omega_j=1$, such that the points in $\mathcal{A}$ are fully covered.   In other word, the optimal location of the polyellipsoid such that the largest sum of the weighted distances from a the demand point to all the translated foci is minimized. This problem naturally extends the widely studied $1$-center problem in case a single focus is considered.

Observe that the demand point $a\in \A$ belongs to the $x$-translation, for some $x\in \R^d$, of a polyellipsoid $\mathbb{P}_{\U,\mathbf{\omega}}(r)$ if $a \in \mathbb{P}_{x+\U,\mathbf{\omega}}(r)$, where $x+\U=\{x+u: u\in \U\}$, i.e., if $\varphi_{\U a}(x) := \sum_{u \in \mathcal{U}} {\bf \omega}_{u} \|a-u-x\| \leq r$.

Thus, the problem can be formulated as:
\begin{equation}
\min_{x \in \R^d} \max_{a \in \mathcal{A}} \varphi_{\U a}(x). \tag{MEP}\label{mep}
\end{equation}
Note that if $x^*$ is an optimal solution of the above problem, it induces a polyellipsoid with foci $\U^* = \{u+x^*: u \in \U\}$, weights $\omega$, and radius $r^*=\max_{a \in \mathcal{A}} \varphi_{\U a}(x)$. Thus, $x^*$ represents the \textit{optimal} translation of the foci in $\U$ to minimally cover the points in $\A$, with respect to the weighted $\omega$-sum of the functions $\varphi_{\U a}(x)$, $\forall a\in \A$. Observe also, that, as usual in minimax problems, one can reformulate the problem by introducing the auxiliary variable $r$ (representing the minimum radius) as follows:
\begin{align}
\min_{x\in \R^d, r \in \R_+} & \;\;r\nonumber\\
\mbox{s.t. } & \dsum_{u \in \mathcal{U}} {\bf \omega}_{u} \|a-u-x\| \leq r, \forall a \in \mathcal{A}. \nonumber
\end{align}
It easily follows that the problem above can be formulated as the minimization of a linear objective function subject only to linear  or norm-type constraints:
\begin{align}
\min & \;\; r \nonumber\\
\mbox{s.t. } & \dsum_{u \in \U} \omega_{u} d_{ua} \leq r, & \forall a\in \A,\nonumber\\
&  \|a-u-x\| \leq d_{ua}, &\forall u \in \U, \; a\in \A, \label{mep1}\tag{MEP$_1$}\\
 & x \in \R^d,& r \in \R,\nonumber\\
 & d_{ua}\in \R^{|\U|\times |\A|}, &\forall u \in \U, a\in \A.\nonumber
 \end{align}
By the comments above, for $\ell_p$-norms, the nonlinear constraints in \eqref{mep1} can be rewritten as a set of SOC constraints. The case in which the norm is a block norm is even simpler since the constraints can be written as linear constraints. Thus, in the two cases the problem is solvable by interior point algorithms using any of the available optimization solvers.

\begin{figure}[h]
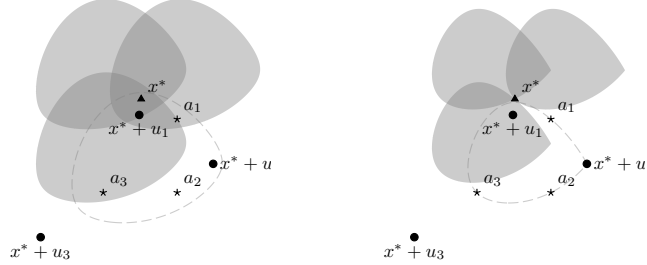

\begin{center}
\input{graph4}~\input{graph5}
\end{center}
\caption{Equivalence on the intersection of the polyellipses centered at the demand points and a feasible covering polyellipse.\label{fig:intersect}}
\end{figure}

Note also that finding a feasible solution, $(x,r)$, to \eqref{mep1} is equivalent to find a point $x$ in the intersection of the $n$ polyellipsoids $\mathbb{P}_{a_1+\U,\mathbf{\omega}}(r), \ldots, \mathbb{P}_{a_n+\U,\mathbf{\omega}}(r)$, since such a solution, $(x,r)$ must verify all the \emph{polyellipsoid} constraints in \eqref{mep1}. In Figure \ref{fig:intersect} (left) we show a toy example with three demand points (stars) to be covered by a (Euclidean) polyellipse with three foci. If we fix a radius and \textit{center} the polyellipse at each of the demand points we obtain the three gray polyellipses in the picture. The point $x^*$ belongs to the intersection of those polyellipses and centering the polyellipse at that point, the polyellipse (dashed line) covers all the demand points (in this case being clearly a non optimal solution). In the right picture we draw the same situation but with the optimal radius of the polyellipse.

The following result allows us to decompose \eqref{mep} into smaller problems:
\begin{thm}\label{theo:1}
Let $r^*$ be the solution of \eqref{mep}. Then, $r^*=\dmin_{S\subset \A:\atop |S|=d+1} r^*_S$, where $r^*_S$ is the optimal solution of the following problem:
\begin{align}
r^*_S := \min & \;\; r  \label{meps} \tag{MEP$_S$}\\
\mbox{s.t. } & \dsum_{u \in \U} \omega_{u} d_{ua} \leq r, & \forall a\in S,\nonumber \\
&  \|a-u-x\| \leq d_{ua}, &\forall u \in \U, \; a\in S,\nonumber\\
 & x \in \R^d, r \in \R,& \nonumber\\
 & d_{ua}\in \R^{|\U|\times |\A|}, &\forall u \in \U, a\in S. \nonumber
\end{align}
for all $S \subseteq \A$ with $|S|=d+1$.
\end{thm}
\begin{proof}
The proof follows by a standard application of Helly's theorem  on the intersection of convex bodies~\cite{helly}.  The optimal solution of \eqref{mep} reduces to finding the smallest $r$ among the solutions of the above collection of subproblems for all the subsets $S \subset \A$ with $|S|=d+1$
\end{proof}

This reformulation of the problem will be used both to derive a polynomial-time complexity for the problem and also to develop an Elzinga-Hearn based solution approach.

\subsection{On the complexity of \eqref{mep}}

In what follows we prove a polynomiality result for the minimum enclosing polyellisoid with a fixed number $k$ of given foci problem in fixed dimension $d$. As usual in computational geometry, the model of computation is that of algebraic computations and comparisons over the reals. The reader is referred to the papers by Renegar~\cite{renegar} and Dyer\cite{Dyer} for further details on the underlying results and the algorithms on which we base our construction.

The main tool to be applied is the linear time algorithm of Dyer for a special class of convex program with few nonlinear constraints. We will prove that an iterative use of that algorithm implies the polynomial-time complexity of the \eqref{mep}. Observe that the reformulation of \eqref{mep} as \eqref{mep1} has a number of nonlinear constraints that are essentially equivalent to polynomials. However, even assuming fixed dimension $d$ and constant number $k$ of foci, the number of nonlinear constraints is \textit{high} as compared with the number of linear ones, and therefore one cannot apply the result of Dyer\cite{Dyer}. In spite of that, we will prove that a convenient decomposition of that problem would allow one to polynomially solve \eqref{mep} in fixed dimension.

\begin{thm}\label{th:1}
Let the number of foci, $k$, and the dimension, $d$, be fixed. Then, \eqref{mep} with Euclidean norm is solvable in polynomial time in $|\A|$.
\end{thm}
\begin{proof}
Recall that by Theorem \ref{theo:1}, solving \eqref{mep1} (and then \eqref{mep}) is equivalent to solve all problems \eqref{meps} for $S \subseteq \A$ with $|S|=d+1$. Observe that in \eqref{meps} the number of nonlinear constraints is constant ($(d+1)\times k$). Let us denote by $K(S)=\{(x,d)\in \mathbb{R}^{d} \times \R^{k \times (d+1)}: \|a-u-x\| \leq d_{ua}, \forall u \in \U, \; a\in S\}$, the feasible domain induced by the nonlinear constraints in \eqref{meps}. Then, we have:
\begin{description}
\item[1. ]  For any $\delta>0$, $a\in S$ and $u\in \U$:
    \begin{equation}
    g_{ua}(x,d):= \|a-u-x\| - d_{ua} \le \varepsilon, \; \forall \varepsilon \in [0,\delta], \mbox{ and } \forall (x,d)\in K(S).\label{item:1}
    \end{equation}

Moreover, since $\|\cdot\|$ is assumed to be the Euclidean norm, the inequality $g_{ua}(x,d)\le \varepsilon$ can be easily transformed into the following polynomial inequality: $ \sum_{j=1}^d (a_j-u_j-x_j)^2 \le (\varepsilon-d_{ua}) ^2$.
\item[2. ] The gradient $\nabla g_{ua}(x,d_{ua})=( - \frac{-a-u+x}{\|a-u-x\|}, -1)^t$. Therefore, the inequality $ \nabla g_{au}^t(x,d_{au}) y\le 0$ for any $y\in \mathbb{R}^{d+1}$ is equivalent to the following polynomial inequality $[\sum_{j=1}^d (a_j-u_j-x_j) y_j]^2 \le \sum_{j=1}^d (a_j-u_j-x_j)^2 y_{d+1}$.
\item[3. ] The largest degree of all the polynomials involved in \eqref{meps} and the transformations made explicitly in items [1.] and [2.] above is 4 (constant) and much smaller than $|\A|$.
\end{description}
Now, we are in position to apply the algorithm of Dyer to problem \eqref{meps} because the number of nonlinear constraints is constant. Thus, that problem is solvable in linear time of the number of linear constrains, $O(d)$, which in this applications turns out to be constant. Hence, because for solving \eqref{mep}, one needs to solve $O(|\A|^{d+1})$ problems in the collection and taking the maximum $r$ value, the original problem can be solved in $O(|\A|^{d+1})$-time.  \hfill $\Box$.
\end{proof}
For the sake of readability, the above result has been presented for the Euclidean norm. However, it extends with minimal changes to any $\ell_p$-norm with $p\in \mathbb{Z}$ and $1<p<+\infty$. Observe that instead of squaring the norms expressions in $g_{ua}$, one can reformulate them as the following set of inequalities:
\begin{align*}
v_{aj} \geq a_j-u_j-x_j, \forall j=1, \ldots, d, a \in S,\\
v_{aj} \geq -a_j+u_j+x_j, \forall j=1, \ldots, d, a \in S,\\
\sum_{j=1}^d v_{aj}^p \le d_{ua}^p, \forall a \in \A.
\end{align*}
where the $v$-variables represent the absolute values of each of the coordinates involved in the norm expressions. Thus, similar arguments to those in the proof of Theorem \ref{th:1} can be derived in case $p$ is small (compared to $|\A|$).

Note also that if $\| \cdot \|$ is polyhedral the constraints in \eqref{mep1} can be written as the following set of linear inequalities:
$$
e^t (a-u-x) \leq d_{ua}, \quad \forall u \in \U, a \in \A, e \in {\rm Ext}_{\|\cdot\|^o}
$$
where ${\rm Ext}_{\|\cdot\|^o} = \{e^o_1, \ldots, e^o_g\}$ are the extreme points of the polar ball of the unit ball of $\| \cdot\|$ (see e.g., \cite{NP06,WW85}). Thus, \eqref{mep} can be solved using linear programming tools.

\section{Solution approaches for solving \eqref{mep}\label{sec:3}}
In this section we describe different solution approaches for solving \eqref{mep} beyond the solution of the SOCP model provided in Section \ref{sec:2}. We describe here two types of approaches. First, we describe those which exploit the strong duality of the problem under two different \textit{formulations}: conic and Lagrangean dual. Then, we provide a geometric construction based on the classical Elzinga-Hearn algorithm for solving the problem.

While the dual approaches allow us to provide alternative mathematical programming formulation for the problem and relate it with other classical problems (as the Weber problem), the Elzinga-Hearn based approach allows, in practice, to efficiently solve geometrically the problem.

\subsection{Conic Dual}

Again, for the sake of simplicity we consider the \eqref{mep} with Euclidean norm. In such a case, the problem can be reformulated as the SOCP \eqref{mep1}. Therefore, one can derive  its conic dual, which results in:
\begin{align}
\max & \sum_{u\in \U} \sum_{a\in A} (a-u) \left(\lambda_{ua}^1 + \lambda_{ua}^2\right) \nonumber\\
\mbox{ s.t. } &  \sum_{a\in A} \mu_a \leq 1,\nonumber\\
& \sum_{u\in \U} \sum_{a\in A} \lambda_{ua}^1 \geq 0,\nonumber\\
&  \sum_{u\in \U} \sum_{a\in A} \lambda_{ua}^2  \leq 0, \label{mep:cdual}\tag{MEP$_{ConicDual}$} \\
& \sqrt{(\lambda_{ua}^1)^2+ (\lambda_{ua}^2)^2}\le \omega_{u} \mu_a, \forall u \in \U,\; a\in A, \nonumber\\
& \mu_a \in \R_+, \lambda^1_{ua}, \lambda^2_{ua} \in \R, \forall a \in \A, u \in \U.\nonumber
\end{align}
Slater condition ensures strong duality between both problems (see e.g. \cite{bv}) giving rise to simple approaches to solve the \eqref{mep}. Again, the case in which the norm is not Euclidean but $\ell_p$ or polyhedral permits a similar representation which results in a SOCP problem (as the above) or a linear problem, respectively.

The main advantage, again, of the primal and the dual SOCP formulations of the problem is that off-the-shelf software packages are capable to solve this type of problems efficiently using interior-point based algorithms. These algorithms applied to SOCP problems are known to have a polynomial time complexity for a given tolerance factor $\varepsilon$, assuring convergence in at most a given number of iterations.

Based on our computational experience, as we will see in Section \ref{sec:5}, the SOCP formulations are not only interesting because its theoretical complexity, but also because current solvers are able to handle medium size instances in a reasonable CPU time with reduced implementation effort.

\subsection{Lagrangean Dual}\label{subsec:lagdual}

Now, we analyze a primal-dual approach based on the Lagrangean dual to solve \eqref{mep}.  For the sake of simplicity, we assume that the norm is strictly convex, although as we will point out at the end of the section, the results extend to non strictly convex norms. Let us denote by $f_a(x) = \sum_{u\in \U} \omega_u \Phi_{ua} (x) = \sum_{u\in \U} \omega_u \|a-u-x\|$, for all $a \in \A$.  Clearly, each $\Phi_{ua}$ is a convex function, and then $f_a$ is convex.

Thus, the Lagrangean dual is:
\begin{align*}
\max_{\alpha \in \R^{|\A|}} \min_{x \in \R^d} & \;\; \dsum_{a \in \mathcal{A}} \alpha_a f_a(x)\\
\mbox{s.t. } &  \dsum_{a \in \mathcal{A}} \alpha_a =1,\\
& \alpha_a \geq 0, \forall a \in \mathcal{A}.
\end{align*}
where $\alpha_a$ for $a\in \A$ are the dual multipliers.

Note that each of the functions $f_a$ attains its minimum in a solution of the Weber problem with demand points $\A_a = \{a-u: u \in \mathcal{U}\}$ and weights $\mathbf{\omega}$, which is unique provided that the norm is strictly convex and the foci in $\mathcal{U}$ are not collinear. Furthermore, if we denote by
$$
F(\alpha) =\min_{x\in \R^d} \dsum_{a \in \mathcal{A}} \alpha_a f_{a}(x) =  \min_{x\in \R^d} \dsum_{a \in \mathcal{A}} \alpha_a \dsum_{u \in \mathcal{U}} {\bf \omega}_{u} \|a-u-x\|,
$$
i.e., the optimal value of the Weber Problem with the $nk$ demand points $\bigcup_{a\in A} \A_a$ and weights $\{\alpha_a\omega_{u}: a\in \mathcal{A}, u \in \mathcal{U}\}$, the dual problem becomes:
\begin{align}
\max_{\alpha \in \R^{|\A|}}  & \;\; F(\alpha)\nonumber\\
\mbox{s.t. } &  \dsum_{a \in \mathcal{A}} \alpha_a =1,\label{dual}\tag{${\rm MEP}_{\rm LagDual}$}\\
& \alpha_a \geq 0, \forall a \in \mathcal{A}.\nonumber
\end{align}
Note that the feasible region is nothing but the $|\mathcal{A}|$-dimensional probabilistic simplex, $\Delta_\A = \{\alpha \in \R_+^{|\mathcal{A}|}:  \dsum_{a \in \mathcal{A}} \alpha_a =1\}$.

By strong duality:
$$
\min_{x \in \R^d} \max_{a \in \mathcal{A}} \dsum_{u \in \mathcal{U}} {\bf \omega}_u \|a-u-x\| = \max_{\alpha \in \Delta_\A} \;\; F(\alpha)
$$
It is straighforward that $F$ is concave on $\Delta_\A$. Observe also that if $\|\cdot\|$ is a strictly convex norm and the points in $\mathcal{A}$ (and in $\mathcal{U}$) are not collinear, for a fixed $\alpha \in \Delta_\A$, the function $h_\alpha(x) = \sum_{a \in \mathcal{A}} \alpha_a f_{a}(x)$, for all $x\in \R^d$, is strictly convex. Denoting by $x_\alpha^* = \arg\min_{x\in \R^d} h_\alpha(x)$ (which is well-defined by the strict convexity of $h_\alpha$) we have that $F(\alpha) = \sum_{a \in \mathcal{A}} \alpha_a f_{a}(x^*_\alpha)$ for all $\alpha \in \Delta_\A$. Furthermore, the gradient of $F$ can be easily derived:
$$
\dfrac{\partial F}{\partial \alpha_a} (\alpha) = f_{a}(x^*_\alpha) =  \dsum_{u \in \mathcal{U}} {\bf \omega}_u \|a-u-x^*_\alpha\|, \quad \forall a \in \mathcal{A}.
$$
With such an information, since \eqref{dual} is an optimization problem with a concave objective function and a linear feasible set, a projected gradient descent method can be applied to solve it. Let start with an initial solution $\alpha^0\in \Delta$ and the iterations are in the form:
$$
\alpha^{k+1} = \Pi_\Delta(\alpha^k - \eta^k \nabla_\alpha F(\alpha^k)).
$$
where $\Pi_\Delta(\beta) = \arg\min_{\gamma \in \Delta} \|\beta-\gamma\|$ is the orthogonal projection onto the unit simplex which can be efficiently computed in $O(n\log{n})$ computation time (see \cite{simplex1,simplex3,simplex2}).

Observe that at each step of this approach one has to compute the gradient $\nabla_\alpha F(\alpha)$ which in turns implies solving the following Weber problem:
$$
\min_{x\in \R^d} \dsum_{a\in \A} \alpha_a \dsum_{u \in \mathcal{U}} {\bf \omega}_u \|a-u-x\|.
 $$
These problems can be arbitrarily approximated using Weiszfeld algorithm and its relatives. Thus, the optimal translation of the covering polyellipsoid, $x$, coincides with an optimal solution of a weighted Weber problem with demand points $\{a-u: a \in \A, u \in \U\}$ and weights $\{\alpha_a\omega_u: a\in \A, u\in \U\}$, and then, as stated in \cite{hansen80}, $x$ belongs to the convex hull of those demand points.

In case the norm $\|\cdot \|$ is not strictly convex, uniqueness of optimal solutions cannot be ensured, but still a similar approach can be followed replacing gradients by subgradients.
Actually,
$$ \partial F(\alpha)= {\rm conv} \left( \bigcup_{\bar{x}\in \arg\min h_{\alpha} (x)}(f_{a}(\bar{x}))_{a\in \A}  \right). $$

\subsection{A decomposition approach for \eqref{mep}}

One of the most popular approaches for solving the planar $1$-center problem is due to Elzinga and Hearn \cite{EH}. The Elzinga-Hearn approach (EH, for short) is based on constructing the minimum enclosing disks of three demand points until the whole sets of points is assured to be fully covered. At each iteration, full coverage is checked, and the triplet of points is changed by incorporating the demand point which is furthest from the actual covering disk replacing one of the three points. The overall worst case complexity of the algorithm is $O(n^3)$, but its performance in practice outperforms other existing strategies. Although the EH algorithm was initially designed to solve the planar unweighted $1$-center problem (taking advantage of the geometry of disks), the approach can be adequately extended to the higher-dimensional weighted case~\cite{HearnVijay}. In what follows we describe how the EH decomposition paradigm can be applied to solve \eqref{mep}.

First, let us assume that $r^*\geq 0$ is the optimal radius for our covering problem \eqref{mep} Since the problem is convex the set of optimal solutions, $X^*$, is convex. Let us denote by:
$$
\P_{\U,\omega}^a(r^*)= \{z \in \R^d: \sum_{u \in \mathcal{U}} \omega_u\|a-u-z\|\leq r^*\}, \forall a \in \mathcal{A},
$$
the polyellipse with foci $\{a-u: u \in \mathcal{U}\}$ and radius $r^*$. Then, any optimal solution $x^*\in \bigcap_{a\in \mathcal{A}} \P_{\U,\omega}^a(r^*)$. Actually, since $r^*$ is minimum:
$$
\bigcap_{a\in \mathcal{A}} \P_{\U,\omega}^a(r^*) = X^*
$$
since otherwise, it would contradict the optimality of $r^*$.

The following result is the main component to apply the decomposition approach to \eqref{mep}.
\begin{thm}\label{thm:3}
There exists  $S \subseteq \mathcal{A}$ with $|S| =d+1$ such that $\displaystyle\bigcap_{a\in S} \P_{\U,\omega}^a(r^*) =  X^*$.
\end{thm}
\begin{proof}
Let $\mathcal{A}_{d+1} = \{S \subseteq \mathcal{A}: |S| =d+1\}$. We define $Q_{r^*}(S)= \bigcap_{a\in S} \P_{\U,\omega}^a(r^*)$ for all $S \in \mathcal{A}_{d+1}$.  Since $\bigcap_{a\in \mathcal{A}} \P_{\U,\omega}^a(r^*) = X^*$, in particular we have that $ X^* \subseteq Q_{r^*}(S) \neq \varnothing$ for all $S \in \mathcal{A}_{d+1}$.

On the other hand, since $\bigcap_{a\in \mathcal{A}}\P_{\U,\omega}^a(r^*)\neq \varnothing$, then, for all sets $S \subseteq \mathcal{A}$ with $|S| =d+1$, $\bigcap_{a\in S} \P_{\U,\omega}^a(r^*) \neq \varnothing$. Suppose that for all $S \in \mathcal{A}_{d+1}$,  $int (Q_{r^*}(S))$  includes $X^*$. Since  $int (Q_{r^*}(S))$ is convex, there exists $\varepsilon>0$ such that  $int(Q_{r^*-\varepsilon}(S)) \neq \varnothing$. Thus, $\bigcap_{a\in S}\P_{\U,\omega}^a(r^*-\varepsilon) \neq \varnothing$, and by Helly's theorem also $\bigcap_{a\in \mathcal{\A}} \P_{\U,\omega}^a(r^*-\varepsilon) \neq \varnothing$, contradicting the optimality of $r^*$.

Alternatively, if there exists $S\in \mathcal{A}_{d+1}$ such that $X^*\subset Q_{r^*}(S)$ and $int (Q_{r^*}(S))=\varnothing$, it implies that $X^*$ is included in one face of dimension at most $d-1$ of $\P_{\U,\omega}^a(r^*)$, for all $a\in S$. Then for any $\varepsilon>0$, $\bigcap_{a\in S}\P_{\U,\omega}^a(r^*-\varepsilon) = \varnothing$, and  this implies that $r^*$ is the minimum value for a full coverage of $\A$ and also that $X^*$ is a set of optimal solutions.
 \hfill $\Box$
\end{proof}

Based on the above result,  a decomposition approach can be derived to solve the problem, as described in Algorithm \ref{alg:EH}.

\begin{algorithm}[h!]
  \renewcommand{\algorithmicrequire}{\textbf{Input:}}
  \renewcommand{\algorithmicensure}{\textbf{Output:}}
  Compute the minimum radius polyellipsoid covering $S^0$: $r^0$, $x^0$.
  \begin{algorithmic}
    \REQUIRE $\mathcal{A}$ and $S^0 \subseteq \mathcal{A}$ with $|S^0|=d+1$. $k=0$.
    \STATE  {\bf 1.} Let $\rho^k = \dmax_{a\in \mathcal{A}} \dsum_{u \in \mathcal{U}} {\bf \omega}_u \|a-u-x^k\|$ and $a^k\in \mathcal{A}$ reaching such a maximum.
    \STATE {\bf 2.} \IF{$\rho^k = r^k$}
    \STATE STOP.
    \ELSE
    \STATE $S^{k+1} = S^k \cup \{a^k\} \backslash \{b^k\}$ with $b^k$ reaching $r^{k+1} = \dmax_{b^k \in S^k} \dmax_{a \in S^k \cup \{a^k\} \backslash \{b^k\}} \dsum_{u \in \mathcal{U}} {\bf \omega}_u \|a-u-x^k\|$
    \ENDIF
    \ENSURE $r^k$ and $x^k$. Go to {\bf 1.}
  \end{algorithmic}
    \caption{A decomposition algorithm for solving \eqref{mep}.  \label{alg:EH}}
\end{algorithm}

\begin{thm}\label{th:coEH}
Algorithm \ref{alg:EH} computes the minimum-radius polyellipsoid in $O(|A|^{d+2})$ time for strictly convex norms.
\end{thm}
\begin{proof}
By Theorem \ref{thm:3}, an optimal polyellipsoid is defined by $d+1$ demand points. Such a polyellipsoid will be determined by the minimum-radius polyellipsoid covering these points. At any iteration $k$ of the algorithm, one computes in {\bf 1.} the minimum-radius polyellipsoid for $d+1$ given points (which takes constant time provided that $d$ is fixed, see Theorem \ref{th:1}). If the optimal polyellipsoid cover all the points, then, the \emph{furthest} points (at distance $\rho^k$) are those at sum of the distances exactly $r^k$. Then, we have constructed a polyellipsoid, that covers all the points in $\mathcal{A}$ and which is the one which minimally covers a subset of $\mathcal{A}$. Otherwise, a new $d+1$-points subset of $\mathcal{A}$ is constructed by replacing one (the most convenient) of its elements by the furthest point to the previous polyellipsoid. Observe that at each iteration the solution of the subproblem is unique because of the strict convexity of the norm. Then, since the number of subsets with $d+1$ elements is $O(|\A|^{d+1})$ and the polyellipsoids generated at each iteration are monotone increasing in radius (the assumption $r^{k+1} > r^k$ is considered because otherwise, the algorithm stops) the algorithm stops in at most $O(|\A|^{d+1})$ iterations and each iteration requires $O(|\A|)$ (the number of comprobations for susbtitution in Step 3). \hfill $\Box$
\end{proof}

Observe that the main advantage of this approach is that one only solves minimum enclosing polyellipsoid problems for $d+1$ points, which in practice is very convenient, in particular when the number of demand points in $\mathcal{A}$ is large. Moreover, in the planar case, no matter the number of demands points in the problem this approach only needs to solve MEP problems covering three of the demand points. Also, we will see that the number of iterations of the approach, as for the $1$-center problem, is small, reducing significatively the computation times as compared with a plain SOCP formulation.

Note also, that in Algorithm \ref{alg:EH}, one has to solve \eqref{mep} for subsets of $\A$ with $d+1$ points. In the planar Euclidean case with a single focus, one can explicilty derive the the minimum enclosing disk of three points~\cite{EH}, being this \textit{oracle} doable in constant time. In the general case, one has to solve a SOC-programming problem, increasing the theoretical complexity of the procedure.

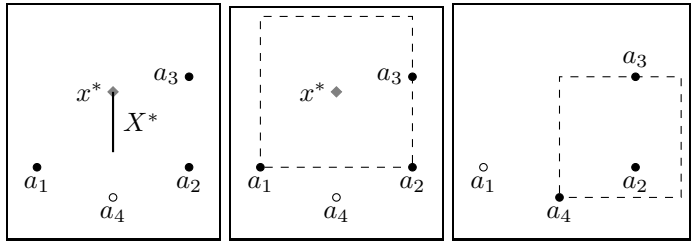
\begin{figure}
\begin{center}
\fbox{\begin{tikzpicture}[scale=2]

\coordinate(a1) at (0,0);
\coordinate(a2) at (1,0);
\coordinate(a3) at (1,0.6);
\coordinate(a4) at (0.5,-0.2);
\coordinate(x) at (0.5,0.5);

\node[circle,draw,fill, inner sep=1](a-1) at (a1) {};
\node[circle,draw,fill, inner sep=1](a-2) at (a2) {};
\node[circle,draw,fill, inner sep=1](a-3) at (a3) {};
\node[circle,draw, inner sep=1](a-4) at (a4) {};
\node[diamond,draw, fill,gray, inner sep=1](x-1) at (x) {};

\node[below] at (a1) {$a_1$};
\node[below] at (a2) {$a_2$};
\node[left] at (a-3) {$a_3$};
\node[below] at (a4) {$a_4$};

\draw[thick] (0.5,0.1)-- node[right] {$X^*$} (0.5,0.5);
\node[left] at (x) {$x^*$};
\node[left] at (1,0.967) {};

\end{tikzpicture}}~\fbox{\begin{tikzpicture}[scale=2]

\coordinate(a1) at (0,0);
\coordinate(a2) at (1,0);
\coordinate(a3) at (1,0.6);
\coordinate(a4) at (0.5,-0.2);
\coordinate(x) at (0.5,0.5);

\node[circle,draw,fill, inner sep=1](a-1) at (a1) {};
\node[circle,draw,fill, inner sep=1](a-2) at (a2) {};
\node[circle,draw,fill, inner sep=1](a-3) at (a3) {};
\node[circle,draw, inner sep=1](a-4) at (a4) {};
\node[diamond,draw, fill,gray, inner sep=1](x-1) at (x) {};

\node[below] at (a1) {$a_1$};
\node[below] at (a2) {$a_2$};
\node[left] at (a-3) {$a_3$};
\node[below] at (a4) {$a_4$};
\node[left] at (x) {$x^*$};

\draw[dashed] (a-1)--(a-2)--(1,1)--(0,1)--(a-1);
\end{tikzpicture}}~\fbox{\begin{tikzpicture}[scale=2]

\coordinate(a1) at (0,0);
\coordinate(a2) at (1,0);
\coordinate(a3) at (1,0.6);
\coordinate(a4) at (0.5,-0.2);

\node[circle,draw,inner sep=1](a-1) at (a1) {};
\node[circle,draw,fill, inner sep=1](a-2) at (a2) {};
\node[circle,draw,fill, inner sep=1](a-3) at (a3) {};
\node[circle,draw,fill, inner sep=1](a-4) at (a4) {};

\node[below] at (a1) {$a_1$};
\node[below] at (a2) {$a_2$};
\node[above] at (a-3) {$a_3$};
\node[below] at (a4) {$a_4$};
\node[left] at (1,0.967) {};
\draw[dashed] (a-4)--(0.5,0.6)--(1.3,0.6)--(1.3,-0.2)--(a-4);
\end{tikzpicture}}
\end{center}
\caption{Example of non-monotonicity in the radius for non strictly convex norms on Algorithm \ref{alg:EH}.\label{fig:monot}}
\end{figure}

In case the norm is not strictly convex, the above procedure is not assured to reach the optimal solution of the problem since the solution of the subproblems with $d+1$ demand points may not be unique. In the example of Figure \ref{fig:monot}, one has four demand points that want to be covered by a square with minimum edge length (i.e., the minimum enclosing polyellipse for $\ell_1$ norm and a single foci). If at the first iteration of Algorithm \ref{alg:EH} one chooses $S^0=\{a_1, a_2, a_3\}$, the optimal centers for the minimal enclosing squares for those three points are drawn in the left picture (solid line). If the solver chooses $x^* \in X^*$ as the optimal center center, clearly, it does not cover $a_4$, so in the next iteration, one element in $S^0$ is replaced by $a_4$. In this case, the three possible replacements result in the same square length. For instance, if $a_1$ is replaced by $a_4$, one gets that one of the minimum enclosing squares for $S^1=\{a_2,a_3,a_4\}$ is the one drawn in the right picture, which has smaller radius than the one obtained in the previous iteration.

To overcome the above undesired fact, we adapt Algorithm \ref{alg:EH}. In this modification, at the $k$th iteration, if $r^{k+1} = \max_{b \in S^k} \dmax_{a \in S^k \cup \{a^k\} \backslash \{b^k\}} \dsum_{u \in \mathcal{U}} {\bf \omega}_u \|a-u-x^k\| < r^k$, then, instead of removing $b^k$ from $S^k$, we keep both $a^k$ and $b^k$ in $S^{k+1}$, increasing by one the cardinality of the set $S^{k}$. Observe that this modification assures that the sequence of radii $\{r^k\}$ is monotone and the algorithm converges in a finite number of steps. Nevertheless, the  complexity of the algorithm is higher now since one may need to perform $O(|\A|^{d+1})$ iterations but in each iteration a  linear problem with up to $O(|\A|)$ constraints could have to be solved.

\section{Computational Experiments\label{sec:5}}

We have run a battery of experiments in order to show the performance of the Second Order Cone programming formulation and the decomposition approach (Algorithm \ref{alg:EH}) for solving \eqref{mep}. We used the classical planar dataset in \cite{eilon50} and also the datasets for location problems recently proposed in \cite{CB2018}. The first instance, consists of the classical $50$-points (\texttt{EWC}) of Eilon, Watson-Gandy and Christofides, while the rest of the instances consist of geographical coordinates of different population areas in Slovakia. The sizes of these instances are $4873$ (\texttt{Partiz\'anske}), 9562 (\texttt{Ko\v{s}ice}), 79612 (\texttt{\v{Z}ilina}) and 663203 (\texttt{Slovakia}). We compute the minimal radius enclosing polyellipse with number of foci ranging in $\{1, 5, 10, 25\}$ and several norms. In particular we use four strictly convex norms: $\ell_{\frac{3}{2}}$, $\ell_2$, $\ell_3$ and $\ell_4$, and three block norms: $\ell_1$, $\ell_\infty$ and \texttt{hex}, where \texttt{hex} is a block norm whose unit ball is a hexagon with extreme points $\{\pm(2, 0), \pm(1, 2), \pm(-1, 2)\}$ (see \cite{NP06}). The foci were randomly chosen from the set of demand points. We analyze both the unweighted and weighted problems. For the latter, we use either random weights (for the EWC instance) or the weights provided in \cite{CB2018} for the rest of the instances.

We implemented both the SOCP model and the decomposition approach in Python 3, using Gurobi 8 as the optimization solver. The experiments were run in a Mac OSX Mojave  with an Intel Core i7 processor at 3.3 GHz and 16GB of RAM.

In Tables \ref{table:1}-\ref{table:3} we report the results of our computational experience. Table \ref{table:1} and Table \ref{table:2} show the CPU times required for solving the problems for $\ell_p$-norm and block norms polyellipses, respectively. We provide, for each instance its size ($|\A|$), the norm used ($\|\cdot\|$) and the number of foci ($|\U|$), the CPU running times (in seconds) for solving \eqref{mep} both with the Second Order Cone programming formulation (\texttt{Time}$_\texttt{SOC}$) and with the decomposition approach (\texttt{Time}$_\texttt{DEC}$). As one can observe, except for the smallest instance ($|\A|=50$),  the decomposition approach outperforms the results obtained with the SOC formulation. Observe that  both approaches optimally find minimum enclosing polyellipsoids but, while the SOC formulation solves the problems handling simultaneously all the points in $\A$ (adding all the constraints enforcing covering all its elements), Algorithm \ref{alg:EH} solves (iteratively) the problems considering  only $d+1$ points at a time (in the strictly convex norm case). This fact is hard to be observed while the number of demand points is small to medium but it clearly pops up as the number of points increases. The reader can also observe in Table \ref{table:3} that the number of iterations needed to solve the problems with the decomposition approach is rather small. This table shows the average number of iterations of the decomposition approach for solving the instances and also the average cardinality of the sets $S^k$ required for solving the block-norm cases. The number of iterations of the decomposition approach for all the instances range in $[2,6]$. For the $\ell_p$-norm polyellipses, it means that at most $6$ \eqref{mep} problems with $3$ demand points had to be solved, to obtain the solution of the original problem. In the case block-norm case, it may happen that one has to increase the number of points for which the \eqref{mep} problem has to be solved (the sets $S^k$). However, the maximum number of points that we obtained in our experiments  was $6$, that combined with the maximum number of iterations, $6$, gives rise to the highly competitive CPU times needed for solving the problems.

We have marked with \texttt{OoM} those instances for which the SOC approach is not able to load/solving the problem because of an \textit{Out of Memory} flag.

We also show the number of iterations performed by our decomposition procedure (\texttt{it}$_\texttt{dec}$), and for the block norms, the maximum cardinality of the sets $S^k$ used in the modification of Algorithm \ref{alg:EH} for non strictly convex norms.

\setlength{\tabcolsep}{1.5pt}
\begin{table}
\begin{center}
\begin{adjustbox}{angle=90}
\begin{tabular}{|c|c|rr|rr|rr|rr||rr|rr|rr|rr|}\cline{3-18}
      \multicolumn{2}{c|}{}       & \multicolumn{8}{|c||}{Unweighted}     & \multicolumn{8}{c|}{Weighted} \\ \cline{3-18}
   \multicolumn{2}{c|}{}    &  \multicolumn{2}{c}{$\ell_{\frac{3}{2}}$}       & \multicolumn{2}{c}{$\ell_{2}$}       &\multicolumn{2}{c}{ $\ell_{3}$  }      & \multicolumn{2}{c||}{$\ell_{4}$}       & \multicolumn{2}{c}{$\ell_{\frac{3}{2}}$}       & \multicolumn{2}{c}{$\ell_{2}$}       &\multicolumn{2}{c}{ $\ell_{3}$  }      & \multicolumn{2}{c|}{$\ell_{4}$}     \\\hline
\multicolumn{1}{|c|}{ $|\A|$ } & \multicolumn{1}{c}{$|\U|$} & \multicolumn{1}{|c}{t$_{SOC}$ } & \multicolumn{1}{c|}{t$_{DEC}$} &\multicolumn{1}{c}{t$_{SOC}$ } & \multicolumn{1}{c|}{t$_{DEC}$} & \multicolumn{1}{c}{t$_{SOC}$ } & \multicolumn{1}{c|}{t$_{DEC}$}  &\multicolumn{1}{c}{t$_{SOC}$ } & \multicolumn{1}{c||}{t$_{DEC}$} &\multicolumn{1}{c}{t$_{SOC}$ } & \multicolumn{1}{c|}{t$_{DEC}$}  &\multicolumn{1}{c}{t$_{SOC}$ } & \multicolumn{1}{c|}{t$_{DEC}$}  &\multicolumn{1}{c}{t$_{SOC}$ } & \multicolumn{1}{c|}{t$_{DEC}$} &\multicolumn{1}{c}{t$_{SOC}$ } & \multicolumn{1}{c|}{t$_{DEC}$} \\\hline
 \multirow{4}{*}{50}    & 1     & 0.01  & 0.14  & 0.00  & 0.02  & 0.01  & 0.01  & 0.01  & 0.02  & 0.01  & 0.02  & 0.00  & 0.02  & 0.01  & 0.01  & 0.01  & 0.02 \\
      & 5     & 0.07  & 0.10  & 0.04  & 0.05  & 0.08  & 0.10  & 0.05  & 0.10  & 0.07  & 0.10  & 0.03  & 0.08  & 0.06  & 0.13  & 0.05  & 0.07 \\
      & 10    & 0.14  & 0.14  & 0.04  & 0.14  & 0.14  & 0.24  & 0.20  & 0.20  & 0.18  & 0.14  & 0.03  & 0.09  & 0.16  & 0.15  & 0.19  & 0.09 \\
      & 25    & 0.56  & 0.36  & 0.12  & 0.20  & 0.57  & 0.49  & 0.48  & 0.47  & 0.53  & 0.55  & 0.14  & 0.43  & 0.44  & 0.35  & 0.57  & 0.64 \\\hline
 \multirow{4}{*}{4873}   & 1     & 1.88  & 0.12  & 0.62  & 0.15  & 2.80  & 0.16  & 3.09  & 0.13  & 1.74  & 0.14  & 0.58  & 0.17  & 2.46  & 0.13  & 2.51  & 0.17 \\
      & 5     & 11.99 & 0.59  & 5.29  & 1.08  & 13.59 & 0.90  & 11.80 & 0.77  & 11.78 & 0.55  & 3.79  & 0.51  & 11.57 & 0.76  & 12.72 & 0.77 \\
      & 10    & 31.61 & 1.79  & 7.28  & 1.05  & 31.53 & 1.60  & 37.75 & 1.64  & 31.87 & 1.50  & 15.20 & 1.07  & 45.39 & 1.49  & 46.34 & 2.03 \\
      & 25    & 117.96 & 2.77  & 23.30 & 3.34  & 112.68 & 3.88  & 129.25 & 3.53  & 148.39 & 2.60  & 39.00 & 2.45  & 167.49 & 3.48  & 147.35 & 3.83 \\\hline
 \multirow{4}{*}{9562}   & 1     & 4.38  & 0.30  & 1.56  & 0.29  & 4.36  & 0.30  & 4.25  & 0.35  & 4.12  & 0.33  & 1.38  & 0.30  & 4.19  & 0.32  & 4.16  & 0.37 \\
      & 5     & 30.58 & 1.27  & 7.11  & 1.38  & 27.46 & 1.70  & 35.00 & 1.36  & 33.42 & 1.28  & 9.90  & 1.30  & 50.72 & 1.33  & 42.58 & 1.38 \\
      & 10    & 97.25 & 2.82  & 16.65 & 1.80  & 71.98 & 2.04  & 84.53 & 2.67  & 205.41 & 3.19  & 33.87 & 2.47  & 245.38 & 2.67  & 160.70 & 3.33 \\
      & 25    & 213.34 & 6.12  & 47.93 & 5.88  & 317.86 & 6.27  & 365.62 & 6.81  & 503.54 & 6.35  & 368.24 & 6.12  & 531.46 & 6.41  & 463.49 & 6.55 \\\hline
 \multirow{4}{*}{79612} & 1     & 52.86 & 1.72  & 16.51 & 1.67  & 56.13 & 1.80  & 67.14 & 1.91  & 52.48 & 1.73  & 16.56 & 1.72  & 55.40 & 1.80  & 66.80 & 1.80 \\
      & 5     & 379.13 & 10.04 & 84.89 & 9.75  & 432.17 & 10.71 & 413.29 & 10.58 & 7200 & 12.77 & 377.48 & 7.72  & 1490.33 & 7.69  & 1376.67 & 7.85 \\
      & 10    & 1427.21 & 29.41 & 186.11 & 19.63 & 1262.71 & 15.21 & 1637.17 & 16.29 & 4448.48 & 29.64 & 1040.49 & 15.45 & 2496.65 & 17.05 & 3479.00 & 17.47 \\
      & 25    & 7200 & 71.20 & 389.17 & 58.66 & 7200 & 38.77 & 7200 & 40.07 & 7200 & 38.43 & 2734.22 & 37.40 & 7200 & 37.55 & 7200& 37.64 \\\hline
   \multirow{4}{*}{663203} & 1 & 1125.80 & 14.39 & {271.17} & 18.82 & {788.36} & 14.70 & {951.78} & 15.22 &  1149.75 & 13.57 & {281.24} & 17.97 & {806.99} & 14.09 & {949.27} & 14.41 \\
         & 5 & {7200.01} & 84.43 & {1785.31} & 64.11 & \texttt{OoM}  & 75.52 & \texttt{OoM}  & 76.07 & \texttt{OoM}  & 80.27 & \texttt{OoM}  & 58.62 & \texttt{OoM}  & 60.80 & \texttt{OoM}  & 62.82 \\
         & 10 & \texttt{OoM}  & 189.41 & \texttt{OoM}  & 231.85 & \texttt{OoM}  & 193.97 & \texttt{OoM}  & 148.38 & \texttt{OoM}  & 157.83 & \texttt{OoM}  & 154.26 & \texttt{OoM}  & 119.90 & \texttt{OoM}  & 123.95 \\
         & 25 & \texttt{OoM}  & 320.64 & \texttt{OoM}  & 434.61 & \texttt{OoM}  & 350.80 & \texttt{OoM}  & 368.30 & \texttt{OoM}  & 292.14 & \texttt{OoM}  & 288.67 & \texttt{OoM}  & 296.56 & \texttt{OoM}  & 305.61\\\hline
\end{tabular}
\end{adjustbox}
\end{center}
\caption{CPU Times for MEP with strictly convex norms.\label{table:1}}
\end{table}
\begin{table}
\begin{center}
\begin{adjustbox}{angle=90}
\begin{tabular}{|c|c|rr|rr|rr||rr|rr|rr|}\cline{3-14}
      \multicolumn{2}{c|}{}       & \multicolumn{6}{|c||}{Unweighted}     & \multicolumn{6}{c|}{Weighted} \\ \cline{3-14}
   \multicolumn{2}{c|}{}    &  \multicolumn{2}{c}{$\ell_{1}$}       & \multicolumn{2}{c}{$\ell_{\infty}$}       &\multicolumn{2}{c||}{ \texttt{hex} }      & \multicolumn{2}{c}{$\ell_{1}$}       & \multicolumn{2}{c}{$\ell_{\infty}$}       &\multicolumn{2}{c|}{ \texttt{hex} }     \\\hline
\multicolumn{1}{|c|}{ $|\A|$ } & \multicolumn{1}{c}{$|\U|$} & \multicolumn{1}{|c}{t$_{SOC}$ } & \multicolumn{1}{c|}{t$_{DEC}$} &\multicolumn{1}{c}{t$_{SOC}$ } & \multicolumn{1}{c|}{t$_{DEC}$} &\multicolumn{1}{c}{t$_{SOC}$ } & \multicolumn{1}{c||}{t$_{DEC}$} &\multicolumn{1}{c}{t$_{SOC}$ } & \multicolumn{1}{c|}{t$_{DEC}$}  &\multicolumn{1}{c}{t$_{SOC}$ } & \multicolumn{1}{c|}{t$_{DEC}$}  &\multicolumn{1}{c}{t$_{SOC}$ } & \multicolumn{1}{c|}{t$_{DEC}$} \\\hline
 \multirow{4}{*}{50}      & 1     & 0.00  & 0.02  & 0.00  & 0.02  & 0.00  & 0.05  & 0.00  & 0.08  & 0.00  & 0.02  & 0.00  & 0.05 \\
      & 5     & 0.01  & 0.03  & 0.01  & 0.05  & 0.02  & 0.17  & 0.01  & 0.03  & 0.01  & 0.03  & 0.03  & 0.09 \\
      & 10    & 0.02  & 0.09  & 0.01  & 0.09  & 0.04  & 0.18  & 0.02  & 0.08  & 0.01  & 0.06  & 0.05  & 0.18 \\
        & 25    & 0.10  & 0.09  & 0.06  & 0.16  & 0.06  & 0.46  & 0.10  & 0.14  & 0.07  & 0.26  & 0.06  & 0.45 \\\hline
 \multirow{4}{*}{4873}  & 1     & 0.17  & 0.16  & 0.08  & 0.19  & 0.21  & 0.52  & 0.18  & 0.14  & 0.10  & 0.18  & 0.22  & 0.51 \\
      & 5     & 0.92  & 0.61  & 0.86  & 0.44  & 1.44  & 3.13  & 0.93  & 0.73  & 0.95  & 0.71  & 1.33  & 2.44 \\
      & 10    & 2.05  & 1.17  & 2.08  & 1.13  & 2.87  & 5.04  & 2.15  & 1.46  & 1.81  & 1.10  & 2.43  & 4.83 \\
      & 25    & 7.06  & 3.61  & 5.64  & 2.79  & 8.63  & 9.10  & 6.50  & 2.11  & 5.59  & 3.08  & 8.78  & 9.06 \\\hline
 \multirow{4}{*}{9562}  & 1     & 0.34  & 0.44  & 0.19  & 0.41  & 0.44  & 0.98  & 0.35  & 0.41  & 0.20  & 0.41  & 0.44  & 1.03 \\
      & 5     & 2.31  & 1.14  & 1.86  & 1.09  & 2.85  & 5.88  & 2.58  & 0.83  & 2.40  & 1.11  & 2.97  & 4.69 \\
      & 10    & 5.62  & 3.34  & 4.17  & 2.10  & 6.45  & 7.11  & 3.01  & 2.78  & 2.83  & 2.76  & 4.21  & 11.69 \\
      & 25    & 15.08 & 6.63  & 11.54 & 5.08  & 20.35 & 35.70 & 13.16 & 6.63  & 10.68 & 6.64  & 17.86 & 29.19 \\\hline
 \multirow{4}{*}{79612} & 1     & 5.25  & 2.65  & 1.55  & 2.64  & 5.66  & 6.06  & 6.75  & 2.59  & 1.74  & 3.21  & 5.56  & 6.00 \\
      & 5     & 26.92 & 11.13 & 22.10 & 8.65  & 34.63 & 47.97 & 29.45 & 8.97  & 18.69 & 6.43  & 32.93 & 28.75 \\
      & 10    & 68.77 & 13.15 & 43.87 & 12.78 & 94.67 & 57.86 & 71.63 & 17.19 & 40.54 & 12.35 & 88.22 & 56.71 \\
      & 25    & 189.67 & 31.70 & 143.21 & 30.64 & 308.12 & 243.43 & 144.43 & 32.56 & 110.38 & 32.37 & 224.39 & 175.40 \\\hline
           \multirow{4}{*}{663203} & 1    & {55.84} & 26.05 & {16.00} & 17.37 & {56.48} & 47.34 & {53.39} & 26.07 & {23.01} & 20.67 & {57.84} & 44.19 \\
         & 5    & {390.16} & 53.92 & {364.48} & 54.29 & {798.43} & 219.50 & {1232.96} & 69.97 & {533.16} & 52.13 & {1234.53} & 213.33 \\
         & 10   & {3674.59} & 141.34 & \texttt{OoM}  & 103.50 & \texttt{OoM}  & 578.83 & \texttt{OoM}  & 135.59 & \texttt{OoM}  & 102.05 & \texttt{OoM}  & 421.48 \\
         & 25   & \texttt{OoM}  & 444.92 & \texttt{OoM}  & 254.45 & \texttt{OoM}  & 1087.21 & \texttt{OoM}  & 332.15 & \texttt{OoM}  & 252.44 & \texttt{OoM}  & 1046.51 \\\hline
\end{tabular}
\end{adjustbox}
\end{center}
\caption{CPU Times for MEP with polyhedral norms.\label{table:2}}
\end{table}

\begin{table}
\begin{center}
\begin{tabular}{|c|c|c|c|c|}
\cline{3-5}
\multicolumn{2}{c}{} & \multicolumn{1}{|c}{$\ell_p$ norms} & \multicolumn{2}{|c|}{block norms} \\\hline
 $|\A|$    &   $|\U|$ &  it$_{DEC}$ &  it$_{DEC}$ & $|S^k|$ \\
\hline
\multirow{4}{*}{50} & 1     & 2.8   & 4.7   & 4.7 \\
      & 5     & 3.8   & 3.5   & 3.0 \\
      & 10    & 3.4   & 3.5   & 3.2 \\
      & 25    & 3.9   & 3.2   & 3.0 \\
\hline
\multirow{4}{*}{4873} & 1     & 3.5   & 4.3   & 4.7 \\
      & 5     & 3.8   & 4.3   & 3.7 \\
      & 10    & 4   & 4.2   & 3.3 \\
      & 25    & 3.6   & 3.7   & 3 \\
\hline
\multirow{4}{*}{9562} & 1     & 4.0   & 5.3   & 4.7 \\
      & 5     & 4.1   & 4.0   & 3.0 \\
      & 10    & 4.0   & 4.7   & 3.7 \\
      & 25    & 4.0   & 5.0   & 3.0 \\
\hline
\multirow{4}{*}{79612} & 1     & 3   & 4.3   & 4.3 \\
     & 5     & 3.8   & 4   & 3.7 \\
     & 10    & 3.9   & 3.2   & 3 \\
     & 25    & 3.6   & 3.5   & 3.2 \\\hline
        \multirow{4}{*}{663203} & 1    & 3.3  & 4.5  & 4.5\\
 & 5    & 3.3  & 3.2  & 3.2\\
& 10   & 4.0  & 3.5  & 3.0\\
 & 25   & 3.3  & 3.5  & 3.0\\
\hline
\end{tabular}
\end{center}
\caption{Average number of iterations in the decomposition approach for all the instances.\label{table:3}}
\end{table}

\section{Case study: One dimensional covering polyellipsoids}\label{sec:4}

In this section we analyze problem \eqref{mep} in case $d=1$, and describe closed formulas for the minimal covering one-dimensional polyellipsoid covering a finite set of demand points. Observe that while the solution of the $1$-center problem on the line can be trivially obtained ($x^*=\frac{1}{2}(\max_{a\in \A} a + \min_{a\in \A} a)$ and $r=\frac{1}{2}(\max_{a\in \A} a - \min_{a\in \A} a)$), a further analysis is needed when the number of foci is greater than $1$.

Let us first analyze the shape of a polyellipsoid on the real line. Let $\U\subset \R$, $\mathbf{\omega} \in \R^{|\U|}_+$ and $r\geq 0$. The one-dimensional polyellipse with foci $\U$, weights $\mathbf{\omega}$ and radius $r$ is:
$$
\mathbb{P}_{\U,\mathbf{\omega}}(r) = \{z\in \R: \sum_{u\in \U} \omega_u |z-u| \leq r\}.
$$
By Proposition \ref{prop:1}, since polyellipsoids are closed and bounded sets, $\mathbb{P}_{\U,\mathbf{\omega}}(r)$ is, either empty, a single point or an interval. On the other hand, the Weber problem on the line, i.e., $\min_{x\in \R} \sum_{u\in \U} \omega_u |x-u|$, is solved at the median interval, $X^*=[x_0^*, x_f^*]$. Let us also denote by $r^* =\sum_{u\in \U} \omega_u |x^*-u|$ for $x^* \in X^*$ and $\bar u = \sum_{u\in \U} \omega_u u$.

The following result whose proof is straightforward provides the explicit shape of $\mathbb{P}_{\U,\mathbf{\omega}}(r)$.
\begin{lem}\label{lemma:2}
There exist  $u^0, u^f \in \U$ such that
$$
\mathbb{P}_{\U,\mathbf{\omega}}(r)=\left\{\begin{array}{cl}
\varnothing & \mbox{$r<r^*$},\\
\left[u^0, u^f\right] & \mbox{ if $r=r^*$,}\\
\left[\dfrac{r - \bar u + 2\sum_{u\in \U_0} \omega_u u}{2\sum_{u\in \U_0} \omega_u- 1},\dfrac{r - \bar u+ 2\sum_{u\in \U_f} \omega_u u}{2\sum_{u\in \U_f} \omega_u -1}\right] & \mbox{if $r>r^*$}.
\end{array}\right.
$$
where $\U_0=\{u\in \U: u\leq u^0\}$ and $\U_f=\{u\in \U: u\leq u^f\}$.
\end{lem}
\begin{proof}
The first statement is straightforward. Let us analyze the case $r\geq r^*$.\\
Assume, w.l.o.g., that the elements in $\U$ are sorted in increasing order, i.e. $\U =\{u_1 < \cdots < u_k \}$, and denote by $u_0 = -\infty$ and $u_{k+1}=\infty$, and $\omega_j = \omega_{u_j}$, for $j=1, \ldots, k$.  For any $z \in \R$, let $s=s(z)\in \{0,1, \ldots, k\}$ such that $u_{s} < z \leq u_{s+1}$, then:
\begin{eqnarray}
\begin{split}
\dsum_{u \in \U} \omega_u |z-u| &= \dsum_{j=1}^{s} \omega_j (z-u_j) -  \dsum_{j=s+1}^{k} \omega_j (z-u_j)\nonumber \\
&= (\dsum_{j=1}^s \omega_j - \dsum_{j=s+1}^k \omega_j) z - \dsum_{j=1}^{s} \omega_j u_j  + \dsum_{j=s+1}^{k}\omega_j u_j\nonumber\\
&= (2\dsum_{j=1}^s \omega_j -1) z  + \bar u -  2\dsum_{j=1}^{s} \omega_j u_j.\label{lemma2:a}
\end{split}
\end{eqnarray}
The extremes of the nonempty closed interval $\mathbb{P}_{\U,\mathbf{\omega}}(r)$  are  those $z \in \R$ verifying that $\dsum_{u\in \U} \omega_u |z-u| =r$. Thus, $z$ is an extreme of $\mathbb{P}_{\U,\mathbf{\omega}}(r)$ if \eqref{lemma2:a} is equal to $r$.
\begin{enumerate}
\item If $\dsum_{j=1}^s \omega_j = \dfrac{1}{2}$ ($= \dsum_{j=s+1}^k \omega_j$), then $r = \dsum_{u\in \U} \omega_u |z-u| = \bar u -  2\dsum_{j=1}^{s} \omega_j u_j$. In case $|\U|$ is even, all the points in $[u_s, u_{s+1}]$ have the same sum of weighted distances to the foci. Otherwise, $\mathbb{P}_{\U,\mathbf{\omega}}(r) = \{u_s\}$. Thus, $u^0=u_s$ and $u^f=u_{s+1}$ or $u^f=u_s$ in the statement of the Lemma. Note that in this case $r$ coincides with $r^*$ and $\mathbb{P}_{\U,\mathbf{\omega}}(r)$ is the polyellipsoid induced by the $\omega$-weighted median of the points.
\item If $\dsum_{j=1}^s \omega_j \neq \dfrac{1}{2}$,  then $z = \dfrac{r - \bar u +  2\dsum_{j=1}^{s} \omega_j u_j}{(2\dsum_{j=1}^s \omega_j -1)}$. In this case, $\mathbb{P}_{\U,\mathbf{\omega}}(r)$ is a proper closed interval with nonempty interior, and then there exist two points defining the extremes of the interval, thus, there exist $s_0$ and $s_f$ in with $s_0 < s_f$ such that, $u^0 = u_{s_0}$ and $u^f = u_{s_f}$ and for which the equation in the Lemma holds.\hfill $\Box$
\end{enumerate}

\end{proof}

Le us consider that a set of demand points $\A \subseteq \R$ is given. Finding the minimum enclosing one-dimensional polyellipsoid with foci $\U \subseteq \R$, consists of finding the minimum value of $r$ such that all the elements in $\A$ belong to a $x$-translation of $\mathbb{P}_{\U,\mathbf{\omega}}(r)$, for some $x \in \R$. Since $r$ determines the length of the interval, it is clear that $\mathbb{P}_{x+\U,\mathbf{\omega}}(r)= [\min_{a\in \A} a, \max_{a\in \A} a]$, since otherwise it would not cover the points or it would not be  minimal.

Applying Lemma \ref{lemma:2} to the $x$-translated polyellipse, we get that, after some algebra, the $x$-translated polyellipsoid $\mathbb{P}_{x+\U,\mathbf{\omega}}(r)$ is in the form:
$$
\left[x + \dfrac{r - \bar u+ 2\dsum_{u\in \U:\atop x+u\leq a^0} \omega_u u }{2\dsum_{u\in \U:\atop x+u\leq a^0} \omega_u- 1},x + \dfrac{r - \bar u+ 2\dsum_{u\in \U:\atop x+u\leq a^f} \omega_u }{2\dsum_{u\in \U:\atop x+u\leq a^f} \omega_u -1}\right]
$$
Equaling the extremes of that interval to $[\min_{a\in \A} a, \max_{a\in \A} a] =: [a^0, a^f]$ we get that:
$$
r^* = (1-2\dsum_{u\in \U:\atop x+u\leq a^0} \omega_u)(a^f-a^0-\dfrac{\dsum_{u\in \U:\atop a^0 <x^*+ u\leq a^f} \omega_u u}{\dsum_{u\in \U:\atop a^0<x^*+u\leq a^f} \omega_u}) + \bar u^\omega- 2\dsum_{u\in \U:\atop x^*+u\leq a^0} \omega_u u - (a^f-a^0)\dfrac{(2\dsum_{u\in \U:\atop x^*+u\leq a^0} \omega_u -1)^2}{2\dsum_{u\in \U:\atop a^0 <x^*+u\leq a^f} \omega_u}
$$
and
$$
x^* = a^0 - \dfrac{r^* - \bar u+ 2\dsum_{u\in \U:\atop x^*+u\leq a^0} \omega_u u}{2\dsum_{u\in \U:\atop x^*+u\leq a^0} \omega_u -1}  = a^f - \dfrac{r^* - \bar u+ 2\dsum_{u\in \U:\atop x^*+u\leq a^f} \omega_u u}{2\dsum_{u\in \U:\atop x^*+u\leq a^f} \omega_u -1}.
$$
Observe also that, although explicit solutions are detailed above for the one-dimensional minimum enclosing polyellipse problem, they depend of the number of foci $x^*+u$ in $(-\infty,a^0]$ and $(-\infty,a^f]$. Thus the procedure to obtain the optimal polyellipse would consist of iterating on the possible values of $s_0 = |\{u\in \U: x^*+u\leq a^0\}|$ and $s_f = |\{u\in \U: x^*+u\leq a^0\}|$, with $s_0, s_f \in \{0, \ldots, |\U|\}$ with $s_0 < s_f$.

However, since $\U$ is a given set of foci, which range, $\max_{u \in \U} u - \min_{u\in \U} u$, is constant, and also $a^f- a^0$, the search can be reduced. In particular, if $\U=\{u_1 < \cdots < u_k\}$, the optimal solution, $x^*$, is valid for the problem, only if:
\begin{align*}
x^*+u_{s_0} \leq a^0 < x^*+u_{s_0+1},\\
x^*+u_{s_f} \leq a^f < x^*+u_{s_f+1}.
\end{align*}
Equivalently, if $x^* \in (a_0-u_{s_0+1},a_0-u_{s_0}] \cap  (a_f-u_{s_f+1},a_0-u_{s_f}]$. Thus, one can restrict the search to those pairs $(s_0,s_f)\in \{0, \ldots, k\}\times \{0, \ldots, k\}$ such that the intersection of the two intervals is nonempty, i.e.,
$$
u_{s_f} - u_{s_0+1} < a^f-a^0 < u_{s_f+1} - u_{s_0}.
$$
By the above relations, given a valid pair $(s_0,s_f)$, one has that $u_{s_0} < u_{s_f}-(a^f-a^0) \leq u_k - (a^f-a^0)$. Thus, one can restrict $s_0$ to $\{0, \ldots, q_0\}$ with $q_0 = \max\{ j\in\{1, \ldots, k\}: u_j \leq u_k-(a^f-a^0)\}$.

\begin{lem}
There exists an unique pair $(s_0,s_f)\in \{0, \ldots, q_0\} \times \{s_0+1, \ldots, k\}$ such that $(x^*,r^*)$ is the optimal center and radius of the minimal enclosing polyellipsoid.
\end{lem}
\begin{proof}
The proof follows by noting that the problem is equivalent to find $x^* \in \mathbb{P}_{a^0 - \U,\mathbf{\omega}}(r)\cap \mathbb{P}_{a^f -\U,\mathbf{\omega}}(r)$. Thus, if two different solutions exist, it contradicts the minimallity of $r$. Thus, $x^*$ is unique, so there exists an unique combination of indices $(s_0,s_f)$ fullfiling the conditions. \hfill $\Box$
\end{proof}

The above result allows us to terminate the search as soon as a solution $x^*$ in the form above is found verifying that $\dmax \{j \in \{1, \ldots, k\}:  u_j \leq a^0 - x\} = s_0$.
\begin{cor}
Provided that the sets $\A$ and $\U$ are sorted, the one-dimensional minimum enclosing polyellipse can be found in $O(|\U|)$ time for any number of points.
\end{cor}

Also, if all the demand points fall inside the interval determined by some translation of the foci, the optimal polyellipsoid can be explicitly stated.
\begin{cor}
If $\max_{u \in \U} u - \min_{u\in \U} u \leq \max_{a \in \A} a- \min_{a\in \A} a$. The minimum enclosing polyellipse is determined by:
$$
r^* = \dfrac{ \max_{a \in \A} a- \min_{a\in \A} a}{2} \mbox{ and }
x^* = \dfrac{ \max_{a \in \A} a + \min_{a\in \A} a}{2} - \dsum_{u\in \U} \omega_u u
$$
\end{cor}

\section{Extensions\label{sec:6}}

This section addresses two extensions of \eqref{mep}   where we can still exploit the methods and tools previously developed. In particular,  we  analyze the problem that, apart from finding the position and the radius of the covering polyellipsoid, incorporates the foci selection among a given set of potential candidates. Also, we extend the notion of polyellipsoid to ordered median polyellipsoid, following the analogy between the Weber and Ordered Median Location problems.
\subsection{Foci Selection}
In many situations, the foci, which are assumed to be fixed in \eqref{mep}, may be unknown. In the following we address the question of how to \emph{optimally} select a given number of foci from a finite set of candidate points.

Let $\mathcal{A} \subseteq \R^d$ be  the set of demand points and $\mathcal{B}\subseteq \R^d$ a set of potential foci. The goal is to select from $\mathcal{B}$  a subset of $k$ foci $\mathcal{U} \subseteq \mathcal{B}$ with $|\mathcal{U}|=k$ to adequately cover the set of demand points by the translated polyellipsoid with foci in $\mathcal{U}$ and minimal radius. The problem can be stated, in a mathematical programming manner, as follows:
\begin{equation}\label{mepfs}
\min_{x \in \R^d,\atop
\U \subseteq \mathcal{B}: |\U|=k} \max_{a \in \mathcal{A}} \varphi_{\U a}(x) \tag{MEP-FS}
\end{equation}
Recall that $\varphi_{\U a} = \dsum_{u \in \U} \omega_u \|x-a-u\|$ for all $\U \subseteq \mathcal{B}$ and $a\in \A$.

Observe that once the foci, $\U$, are selected, the problem reduces to \eqref{mep}. However, the enumeration of the ${|\mathcal{B}|} \choose {k}$ combinations of possible foci becomes intractable in practice. Thus, we propose, first, a mixed integer non linear programming formulation for the problem, and then, we develop a decomposition approach for solving the problem geometrically.

In the mathematical programming formulation, apart from the variables used to reformulate \eqref{mep}, we use the following set of binary decision variables:
$$
y_u = \left\{\begin{array}{cl} 1 & \mbox{if $u$ is selected as a focus}\\
0 & \mbox{otherwise.}
\end{array}\right., \quad \forall u \in \mathcal{B}.
$$

The following model,  (${\rm MEPFS}$), allows us to determine the optimal foci and the polyellipsoid.
\begin{align}
\min & \; r\label{prowf:of}\\
\mbox{s.t. } & \dsum_{u \in \mathcal{B}} y_u=k,\label{prowf:c3}\\
& d_{au} \geq \|a - u -x \| y_u, \forall a \in \mathcal{A}, u \in \mathcal{B},\label{prowf:c2}\\
&r \geq \dsum_{u \in \mathcal{B}} \omega_{ua} d_{ua}, \forall a \in \mathcal{A},\label{prowf:c1}\\
& r, d_{ua} \geq 0, x\in \R^d,\nonumber\\
& y_u \in \{0,1\}.\nonumber
\end{align}
The objective function \eqref{prowf:of}, as in \eqref{mep},  minimizes the radius, $r$, by choosing an adequate subset of $k$ foci from $\mathcal{B}$ \eqref{prowf:c3}. The way the distances are accounted in the problem depends on the selected foci: in constraint \eqref{prowf:c2}, the distance between a demand point $a\in \A$ and a translated foci $u+x$ with $u \in \mathcal{B}$ is either $\|a-u-x\|$ if $u$ is chosen from $\mathcal{B}$ (i.e., if $y_u=1$) or $0$ otherwise. Finally, as in \eqref{mep}, constraint \eqref{prowf:c1} ensures that $r$ is defined as the maximum among all the sum of the distances from each demand point to the translated selected foci.

The formulation above corresponds to a mixed integer non-linear program. The discrete character comes from the $y$-variables, whereas the non-linearity appears by the constraints \eqref{prowf:c2}, which can be rewritten using big-$M$ constants as follows:
$$
d_{au} + M_{au} ( 1- y_u) \geq \|a - u -x \|, \forall a \in \mathcal{A}, u \in \mathcal{B},
$$
where $M_{au}$ is an upper bound on the value of the norm $\|a-u-x\|$. As already mentioned in Section \ref{subsec:lagdual}, the optimal $x$-value coincides with the solution of a weighted Weber problem for the  set of the demand points $\{a-u: u \in \U\}$ for some $a\in \A$, thus, these big $M$-constants can be derived explicitly.

The above formulation reduces to a Mixed Integer Second Order Cone Optimization (MISOCO) problem provided that the considered norm is polyhedral or in the  $\ell_p$ ($p\ge 1$) family. Then, medium size instances can be solved with nowadays available off-the-shell software.

In what follows, we describe an adaptation of the Elzinga-Hearn based approach described in Algorithm \ref{alg:EH} to this new case where the foci have to be selected from $\mathcal{B}$. Recall that the success of the EH-based approach comes from decomposing the original problem into smaller ones, by solving the minimal covering polyellipsoid problem on a reduced subset of demand points, in particular on subsets of demand points $S \subseteq \A$ with cardinality $d+1$. Following a similar scheme, for each subset $S \subset \A$ with $|S|=d+1$, one has to compute not only the polyellipsoid but also the optimal $k$ foci from $\mathcal{B}$. Let us denote by $\U_R(S)$ the optimal set of foci from $\mathcal{B}$ that minimally cover  the points in $S$, with the condition that the points in $R \subseteq \mathcal{B}$ are not in $\U_R(S)$, i.e., the solution of the following problem:
\begin{equation}\label{mepfs:S}
r_{S,R} := \min_{\U \subseteq \mathcal{B}\backslash \{R\}:
|\U|=k}   \max_{a \in S} \varphi_{\U a}(x)
\end{equation}
Observe that the problem above is a reduced version of \eqref{mepfs} in which the set of potential foci is $\mathcal{B} \backslash \{R\}$ and the set of demand points is $S$. Our decomposition approach for solving \eqref{mepfs} consists of solving problems in the form of \eqref{mepfs:S} , sequentially for different sets $S$, until a termination covering criterion is met. The pseudocode of this method is described in Algorithm \ref{alg:EH2}.

First, we initialize the set $S$ to a set of $d+1$ demands points from $\A$ and $R=\varnothing$. At each iteration, $\texttt{it}$,  we compute the set of optimal foci (solving \eqref{mepfs:S}) with demand points in the actual set, $S^\texttt{it}$, excluding those in the set $R$. Apart from the set of optimal foci for those points, the problem gives us a lower bound of the optimal radius for the covering polyellipsoid, $r_{S^\texttt{it},R}$, which is updated. Since the original problem is always feasible, the restricted problem \eqref{mepfs:S} is feasible if and only if one can choose $k$ elements from $\mathcal{B} \backslash R$, i.e., if $|\mathcal{B} \backslash R| \geq k$. If this were not possible, the process is finished. Once a set of \textit{optimal} foci, $\U^{\texttt{it}+1}$, is computed, one has to  solve \eqref{mep} for those foci using  Algorithm \ref{alg:EH} which also gives the $d+1$ points defining the covering polyellipsoid (the set $S^{\texttt{it}+1}$ in the pseudocode). The solution of that problem provides an upper bound on the optimal value of \eqref{mepfs}  ($UB_\texttt{it}$) which is updated at each iteration. The procedure is repeated until the upper and lower bounds coincide or the number of  reduced potential foci is less than $k$.
\begin{algorithm}[h!]
  \renewcommand{\algorithmicrequire}{\textbf{Input:}}
  \renewcommand{\algorithmicensure}{\textbf{Output:}}
  \begin{algorithmic}
    \REQUIRE

    \STATE $\mathcal{A}$, $S^0 \subseteq \mathcal{A}$ with $|S^0|=d+1$, $\texttt{it}=0$
    \STATE $R  =\varnothing$, $UB=\infty$, $LB=0$.
     \WHILE{$UB > LB$ and $|\mathcal{B} \backslash R| \geq k$}
     \STATE $\U^{\texttt{it}+1} =\U_{R^{\rm it}}(S^{\rm it})$,  $LB=\max\{LB,r_{S^{\rm it}, R}\}$.
     \STATE Solve \eqref{mep} for $\U=\U^{\texttt{it}+1}$ and set $S^{\texttt{it}+1}$ and  $UB=\min\{UB,UB_\texttt{it}\}$).
     \STATE $R = R \cup \; \; \U^{\texttt{it}+1}$.
     \STATE ${\rm it} \rightarrow \texttt{it}+1$.
     \ENDWHILE
    \ENSURE $\U =  \U^{\texttt{it}}$ and $r=UB$.
  \end{algorithmic}
    \caption{Decomposition approach for solving \eqref{mepfs}. \label{alg:EH2}}
\end{algorithm}
Observe that Algorithm \ref{alg:EH2} terminates in a finite number of iterations since $\mathcal{B}$ and $\A$ are finite sets. Note also that two situations may lead to the termination of the algorithm. On the one hand, if at some iteration $UB \leq LB$, then, it would imply that a feasible solution of \eqref{mep} has been found (with radius $UB$) and that for a subset $S^{\texttt{it}+1} \subseteq \mathcal{A}$, any choice among the non explored potential foci, gets a covering polyellipsoid with a larger radius than one already computed. Thus, the best solution found at that iteration must be optimal. On the other hand, if at some iteration $|\mathcal{B} \backslash R| < k$, then all the possible foci have already been explored, and the solution must be among those already computed. Thus, the algorithm outputs the optimal solution of \eqref{mepfs}, as desired.

To conclude, we would like to add some comments on the complexity of Algorithm \ref{alg:EH2}. At each iteration, for solving \eqref{mepfs}, one has to solve a MISOCO which in general could be NP-hard. Thus, although \eqref{mep} can be solved in $O(|\A|^{d+1})$-time, the overall complexity of Algorithm \ref{alg:EH2} is in general non polynomial.

We have run some experiments in order to analyze the performance of Algorithm \ref{alg:EH2} compared to the MISOCO formulation. We have used two of the datasets used in Section \ref{sec:5}, the one with $50$ points from \cite{eilon50} and the one with $4873$ points from \cite{CB2018}. We use unweighted instances and norms in $\{\ell_1, \ell_{\frac{3}{2}}, \ell_2, \ell_3\}$. For each dataset we have randomly chosen the set of potential foci, $\mathcal{B}$, from the demand points with sizes $10$, $15$ and $20$. We consider the foci to choose, $k$, ranging in $\{5,10\}$, and we set a time limit of $30$ minutes in Gurobi for solving the MISOCO problems.

We report in Table \ref{t:MEPSF} the results of our experiments. We show, for each instance the CPU times devoted to the MISOCO formulation (t$_{\mbox{\tiny MISOCO}}$) and for Algorithm \ref{alg:EH2} (t$_{\mbox{\tiny DEC}}$). For the latter, we also provide the number of iterations needed to solve the problems, i.e. the number of MISOCO problems with three demand points (\texttt{It}) as well as the average number of iterations of Algorithm \ref{alg:EH} each time it is called in the procedure  ( \texttt{it}$_{\mbox{\tiny DEC}}$). W also report the MINLP gap obtained with the MISOCO formulation at the time limit  (\%\texttt{Gap}$_{MISOCO}$).

As one can observe from the results, Algorithm \ref{alg:EH2} outperforms the MISOCO formulation in all the instances. Moreover, in most of the instances for the $4873$-dataset, the solver was not able to find the optimal solution within the time limit with the MISOCO formulation. We observe that loading the MISOCO problem in Gurobi is highly time consuming, particularly when non-Euclidean or polyhedral norms are used because the large number of auxiliary variables, and SOC and linear constraints that have to be introduced to represent the norms. This can be seen in columns  \#\texttt{LinCtrs},  \#\texttt{SOCCtrs} and \# \texttt{BinVars}, that  indicate the number of linear, SOC constraints, and binary variables of the entire MISOCO problem, respectively. We have highlighted, with $\texttt{TL}^*$, the instances for which the MISOCO problem was not able to be loaded in Gurobi within 4 hours  (in those instances the gap is not available).
\begin{table}
  \centering
\begin{tabular}{|c|c|c|c|r|r|c|c|c|c|c|c|}\hline
$|\A|$ & $|\mathcal{B}|$ & $k$ & \texttt{Norm} & t$_{\mbox{\tiny MISOCO}}$ &  t$_{\mbox{\tiny DEC}}$ & \texttt{it} & \texttt{it}$_{\mbox{\tiny DEC}}$ & $\%$\texttt{Gap} &  \#\texttt{LinCtrs} & \#\texttt{SOCCtrs} & \# \texttt{BinVars} \\\hline
\hline
 \multirow{20}{*}{50} &  \multirow{4}{*}{10} & \multirow{8}{*}{5} & $\ell_1$  & 1.21  & 0.22  & 1     & 3     & {0.00\%} & 2500  & 0     &  \multirow{4}{*}{510} \\
\cline{4-11}      &       &       & $\ell_{\frac{3}{2}}$  & 22.46 & 2.11  & 3     & 4     & {0.00\%} & 2500  & 2000  &  \\
\cline{4-11}      &       &       & $\ell_{2}$  & 5.97  & 3.38  & 12    & 3     & {0.00\%} & 2000  & 500   &  \\
\cline{4-11}      &       &       & $\ell_{3}$  & 25.11 & 9.88  & 18    & 4     & {0.00\%} & 2500  & 2000  &  \\
\cline{2-2}\cline{4-12}      & \multirow{8}{*}{15} &       & $\ell_1$  & 5.16  & 1.35  & 1     & 5     & {0.00\%} & 3750  & 0     & \multirow{8}{*}{765} \\
\cline{4-11}      &       &       & $\ell_{\frac{3}{2}}$  & 404.10 & 12.95 & 2     & 3     & {0.00\%} & 3750  & 3000  &  \\
\cline{4-11}      &       &       & $\ell_{2}$  & 109.46 & 3.61  & 3     & 3     & {0.00\%} & 3000  & 750   &  \\
\cline{4-11}      &       &       & $\ell_{3}$  & 319.86 & 32.12 & 13    & 3     & {0.00\%} & 3750  & 3000  &  \\
\cline{3-11}      &       &  \multirow{4}{*}{10} & $\ell_1$  & 3.81  & 0.70  & 1     & 4     & {0.00\%} & 3750  & 0     &  \\
\cline{4-11}      &       &       & $\ell_{\frac{3}{2}}$  & 939.28 & 6.43  & 2     & 4     & {0.00\%} & 3750  & 3000  &  \\
\cline{4-11}      &       &       & $\ell_{2}$  & 208.36 & 2.65  & 2     & 4     & {0.00\%} & 3000  & 750   &  \\
\cline{4-11}      &       &       & $\ell_{3}$  & 516.76 & 121.49 & 26    & 3     & {0.00\%} & 3750  & 3000  &  \\
\cline{2-12}      & \multirow{8}{*}{20} &  \multirow{4}{*}{5} & $\ell_1$  & 13.57 & 3.08  & 1     & 4     & {0.00\%} & 5000  & 0     & \multirow{8}{*}{1020} \\
\cline{4-11}      &       &       & $\ell_{\frac{3}{2}}$  & \texttt{TL}\;\;   & 14.88 & 2     & 4     & {0.00\%} & 5000  & 4000  &  \\
\cline{4-11}      &       &       & $\ell_{2}$  & 517.75 & 4.43  & 2     & 2     & {0.00\%} & 4000  & 1000  &  \\
\cline{4-11}      &       &       & $\ell_{3}$  & \texttt{TL}\;\;   & 24.40 & 3     & 4     & {0.00\%} & 5000  & 4000  &  \\
\cline{3-11}      &       &  \multirow{4}{*}{10} & $\ell_1$  & 211.46 & 3.48  & 1     & 3     & {0.00\%} & 5000  & 0     &  \\
\cline{4-11}      &       &       & $\ell_{\frac{3}{2}}$  & \texttt{TL}\;\;   & 377.09 & 3     & 3     & {0.00\%} & 5000  & 4000  &  \\
\cline{4-11}      &       &       & $\ell_{2}$  & \texttt{TL}\;\;   & 58.46 & 2     & 4     & {0.00\%} & 4000  & 1000  &  \\
\cline{4-11}      &       &       & $\ell_{3}$  & \texttt{TL}\;\;   & 175.92 & 2     & 3     & {0.53\%} & 5000  & 4000  &  \\
\hline
 \multirow{20}{*}{4873} &  \multirow{4}{*}{{10}} & \multirow{8}{*}{5} & $\ell_1$  & 479.20 & 1.77  & 1     & 5     & {0.00\%} & 243650 & 0     &  \multirow{4}{*}{48740} \\
\cline{4-11}      &       &       & $\ell_{\frac{3}{2}}$  & \texttt{TL}\;\;   & 8.95  & 8     & 3     & {0.04\%} & 243650 & 194920 &  \\
\cline{4-11}      &       &       & $\ell_{2}$  & 1465.14 & 1.82  & 2     & 4     & {0.00\%} & 194920 & 48730 &  \\
\cline{4-11}      &       &       & $\ell_{3}$  & \texttt{TL}\;\;   & 2.27  & 2     & 4     & {0.01\%} & 243650 & 194920 &  \\
\cline{2-2}\cline{4-12}      & \multirow{8}{*}{{15}} &       & $\ell_1$  & \texttt{TL}\;\;   & 2.54  & 1     & 6     & {0.00\%} & 365475 & 0     & \multirow{8}{*}{73110} \\
\cline{4-11}      &       &       & $\ell_{\frac{3}{2}}$  & \texttt{TL}\;\;   & 15.20 & 5     & 3     & {0.56\%} & 365475 & 292380 &  \\
\cline{4-11}      &       &       & $\ell_{2}$  & \texttt{TL}\;\;   & 2.58  & 2     & 4     & {0.25\%} & 292380 & 73095 &  \\
\cline{4-11}      &       &       & $\ell_{3}$  & \texttt{TL}\;\;   & 6.13  & 2     & 4     & {0.16\%} & 365475 & 292380 &  \\
\cline{3-11}      &       &  \multirow{4}{*}{10} & $\ell_1$  & \texttt{TL}\;\;   & 2.29  & 1     & 3     & {0.00\%} & 365475 & 0     &  \\
\cline{4-11}      &       &       & $\ell_{\frac{3}{2}}$  & \texttt{TL}$^*$  & 14.10 & 3     & 3     & -     & 365475 & 292380 &  \\
\cline{4-11}      &       &       & $\ell_{2}$  & \texttt{TL}$^*$  & 4.67  & 2     & 4     & -     & 292380 & 73095 &  \\
\cline{4-11}      &       &       & $\ell_{3}$  & \texttt{TL}$^*$  & 12.51 & 3     & 4     & -     & 365475 & 292380 &  \\
\cline{2-12}      & \multirow{8}{*}{{20}} &  \multirow{4}{*}{5} & $\ell_1$  & \texttt{TL}\;\;  & 3.36  & 1     & 4     & 9.19\%     & 487300 & 0     & \multirow{8}{*}{97480} \\
\cline{4-11}      &       &       & $\ell_{\frac{3}{2}}$  & \texttt{TL}$^*$  & 33.81 & 3     & 3     & -     & 487300 & 389840 &  \\
\cline{4-11}      &       &       & $\ell_{2}$  & \texttt{TL}\;\;   & 9.58 & 3     & 3     & 0.02\%    & 389840 & 97460 &  \\
\cline{4-11}      &       &       & $\ell_{3}$  & \texttt{TL}$^*$  & 18.21 & 2     & 4     & -     & 487300 & 389840 &  \\
\cline{3-11}      &       &  \multirow{4}{*}{10} & $\ell_1$  & \texttt{TL}\;\; & 6.21  & 1     & 6     & {0.00\%}    & 487300 & 0     &  \\
\cline{4-11}      &       &       & $\ell_{\frac{3}{2}}$  & \texttt{TL}$^*$  & 115.45 & 3     & 3     & -     & 487300 & 389840 &  \\
\cline{4-11}      &       &       & $\ell_{2}$  & \texttt{TL}\;\;  & 87.89 & 3     & 4     & 2.92\%     & 389840 & 97460 &  \\
\cline{4-11}      &       &       & $\ell_{3}$  & \texttt{TL}$^*$  & 201.74 & 4     & 4     & -     & 487300 & 389840 &  \\
\hline
\end{tabular}
  \caption{Results of the Experiments for the Foci Selection Problem.\label{t:MEPSF}}
\end{table}

\subsection{Ordered Median Polyellipsoids}
As mentioned in the introduction, polyellipsoids are identified with the levels curves of the classical Weber (median) problem. Several extensions of this problem has been analyzed in the Location Analysis literature. One of them is the so-called ordered median continuous location problem~\cite{BPE14,ChiaEtal,NP06}. In this problem, apart from the demand points, $\U = \{u_1, \ldots, u_k\}\in \R^d$, a distance measure induced by a norm $\|\cdot\|$ and weights $\omega_1, \ldots, \omega_k \in \R_+$, one is also given a set of modelling weights $\lambda_1, \ldots, \lambda_k \in \R$ and the goal is to find the optimal placement $x$ minimizing the following aggregating function of the distances:
$$
F_{\U,\mathbf{\omega}}^{\mathbf{\lambda}} (x) = \dsum_{i=1}^k \lambda_{i} \omega_{\sigma(i)} \|x- u_{\sigma(i)}\|,
$$
where $\sigma$ is a permutation of the indices $\{1, \ldots, k\}$ such that $\omega_{\sigma(i)} \|x- u_{\sigma(i)}\|  \geq \omega_{\sigma(i+1)} \|x- u_{\sigma(i+1)}\|$ for all $i=1, \ldots, k-1$. The \emph{ordered median function} generalizes most of the objective functions considered in facility location. For instance, choosing $\lambda=(1, \ldots, 1)$, one gets the classical sum of the distances aggregation criterion (Weber problem), while if $\lambda=(1, 0, \ldots, 0)$ one gets the $\max$ criterion which is used in the center problem. It is known that $F_{\U,\mathbf{\omega}}^{\mathbf{\lambda}}$ is convex if and only if the $\lambda$-weights verify $\lambda_1 \geq \cdots \geq \lambda_{k}$. Thus, we will assume such a condition on the weights.

In what follows we draw some comments on the extension of the notion of polyellipsoid to the case in which they are defined as the level curves of functions in the form of $F_{\U,\mathbf{\omega}}^{\mathbf{\lambda}}$.

\begin{figure}
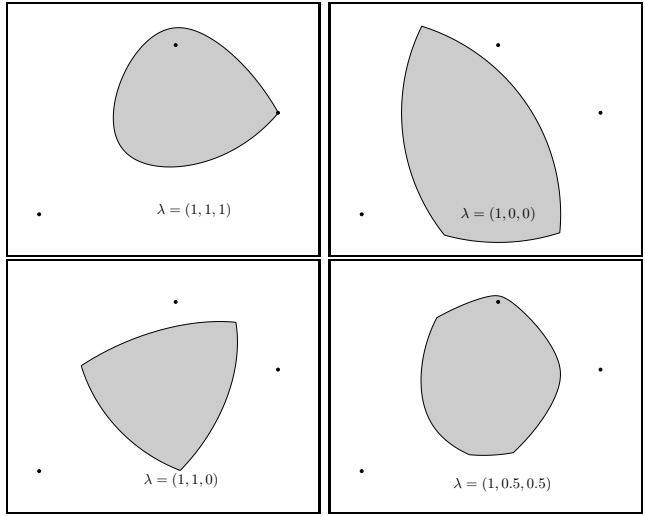

\begin{center}
\fbox{\input{graph6}}~\fbox{\input{graph7}}\\
\fbox{\input{graph8}}~\fbox{\input{graph9}}
\end{center}
\caption{Shapes of ordered median polyellipses for the same foci $\U=\{(-1,5), (2, 3), (-5,0)\}$ and different $\mathbf{\lambda}$-weights.\label{fig:om}}
\end{figure}

Let $\mathcal{U} = \{u_1, \ldots, u_k\} \subseteq \R^d$ be a set of foci and $\omega_1, \ldots, \omega_k \in \R_+$ a set of weights for the foci. Let us also consider a set of weights $\lambda_1, \ldots, \lambda_k \in \R_+$ such that $\lambda_1 \geq \cdots \geq \lambda_{k}$, the $\mathbf{\lambda}$-ordered median polyellipsoid with radius $r$ is defined as:
$$
\E_{\U,\mathbf{\omega}}^\mathbf{\lambda}(r) = \{z \in \R^d: \dsum_{i=1}^k \lambda_{i} \omega_{\sigma(i)} \|z- u_{\sigma(i)}\| = r\}
$$
We denote by $\P_{\U,\mathbf{\omega}}^\mathbf{\lambda}(r)$ the region bounded by $\E_{\U,\mathbf{\omega}}^\mathbf{\lambda}(r)$. Since $F_{\U,\mathbf{\omega}}^{\mathbf{\lambda}}$ is convex, then $\P_{\U,\mathbf{\omega}}^\mathbf{\lambda}(r)$ is a nonempty bounded and convex set for $r\geq \min_{x\in \R^d} F_{\U,\mathbf{\omega}}^{\mathbf{\lambda}} (x)$. In Figure \ref{fig:om} we show four different choices for the lambda weights and the resulting ordered median polyellipses for the same set of three foci on the plane.

Analogously to \eqref{mep}, the ordered median minimal enclosing polyellipsoid problem consists of finding the minimum radius, $r$, and the translation of the ordered median polyellipsoid $\P_{\U,\mathbf{\omega}}^\mathbf{\lambda}(r)$ such that all the demand points are covered by the obtained polyellipsoid. A mathematical programming formulation for the problem is the following:
\begin{equation}\label{mepom}
\min_{x\in \R^d, r \in \R_+} \max_{a\in \A} \dsum_{i=1}^k \lambda_i {\bf \omega}_{\sigma_a(i)} \|a-u_{\sigma_a(i)}-x\|\tag{${\rm MEP}^\lambda_{\rm OM}$}
\end{equation}
where for all $a\in \A$, $\sigma_a$ is a permutation of the indices $\{1, \ldots, k\}$ such that $\omega_{\sigma(i)} \|a-u_{\sigma(i)}^a-x\| \geq  \omega_{\sigma(i+1)}\|a-u_{\sigma(i+1)}^a-x\|$ for all $i=1, \ldots, k-1$.
The problem can be reformulated, using the results proved in \cite{BPE14} as:
\begin{align*}
\min & \;\; r\\
\mbox{s.t. }
& r \geq \dsum_{i=1}^k (v_{ai} + t_{ai}) , \forall a \in \A, \\
& v_{ai} + t_{aj} \geq \lambda_j {\bf \omega}_{i} \|a-u_{i}-x\|, \forall a \in \A, i, j=1, \ldots, k,\\
& r \geq 0,\\
& v_{ai}, t_{ai} \in \R, \forall a \in \A, i=1, \ldots, k.
\end{align*}
where the auxiliary variables $v$ and $t$ allow one to represent the \emph{sorting} factor in the problem without using binary variables. As in \eqref{mep}, in case the norm  is polyhedral or in the $\ell_p$ family with $p\geq 1$, the problem above can be reformulated as a second order cone programming problem.

A further geometrical analysis for this problem is in order. For the sake of presentation, let us assume that the considered norm is strictly convex. The optimal solution of the problem (\ref{mepom}) is given by the smallest $r$ such that $\cap_{a\in \A} \P_{\U-a,\mathbf{\omega}}^\mathbf{\lambda}(r)\neq \varnothing$. Nevertheless, the geometrical structure of the sets $\P_{\U-a,\mathbf{\omega}}^\mathbf{\lambda}(r)$ is intricate since the evaluation depends on the relative sorting of the points with respect to the translated foci. In order to apply an ``a la'' Elzinga-Hearn algorithm to solve this problem we have to introduce the concept of ordered region (see \cite{pufe00}). An ordered region $O_\sigma$ is a closed, connected region of $\mathbb{R}^d$ such that for all $z\in O_\sigma$ a permutation that sorts the vector $(\omega_1\|a-u_1-z\|,\ldots,\omega_k\|a-u_k-z\|)$ in non-increasing sequence is $\sigma$. That is, the chain of inequalities $\omega_{\sigma_1}\|a-u_{\sigma_1}-z\|\ge \omega_{\sigma_2}\|a-u_{\sigma_2}-z\|\ge \ldots \ge \omega_{\sigma_k}\|a-u_{\sigma_k}-z\|$ is valid for all $z\in O_\sigma$. It is straightforward to observe that a subdivision, $\mathcal{S}_a$, of $\mathbb{R}^d$ into ordered regions with respect to the point $a\in \A$ can be obtained with the following arrangement of functions: $\omega_i\|a-u_i-z\|=\omega_j\|a-u_j-z\|, \; \forall i< j\in \{1,\ldots,k\}$.
Denote by $\mathcal{S}$ the intersection of all the subdivisions $\mathcal{S}_a,\; a\in \A$. In each cell of the subdivision $\mathcal{S}$ the order of the vectors $(\omega_1\|a-u_1-z\|,\ldots, \omega_k\|a-u_k-z\|),\; \forall a\in A$ is constant and thus, $\P_{\U-a,\mathbf{\omega}}^\mathbf{\lambda}(r)$ for all $a\in \A$ are standard polyellipsoids. In order to have a valid application of a decomposotion algorithm for problem (\ref{mepom}), similar to Algorithm \ref{alg:EH}, one would have to replicate it on each cell of the subdivision $\mathcal{S}$. This would imply to solve constrained (to the cell) \eqref{mep}  until an optimal restricted solution on that cell is found. Then, choosing the minimum value among all for the different cells will result in the optimal solution for (\ref{mepom}).

To have an idea of the complexity of these subdivisions we present the analysis for the intuitive, important case of the unweighted Euclidean norm. In this case, for each $a\in \A$ the elements that induce the arrangement $\mathcal{S}_a$ are $\|a-u_i-z\|_2=\|a-u_j-z\|_2, \; \forall i< j\in \{1,\ldots,k\}$. It is well-known that for each $i\neq j$ these equations define Euclidean bisectors (hyperplanes) and therefore the number of ordered regions is the number of cells in an arrangement of $O(|\A|^2)$ hyperplanes in dimension $d$. This number is known to be $O(|\A|^{2d})$. Intersecting these subdivision for all $a\in \A$ gives rise to the following upper bound, $O((|\A|^3)^d)$, on the number of elements in the subdivision $\mathcal{S}$. The discussion above allows us to state the following results, whose proof follows from the application of Theorem \ref{th:coEH}, with complexity $O(|\A|^{d+2})$, on each one of the $O(|\A |^{3d})$ full dimensional cells of the subdivision $\mathcal{S}$.
\begin{thm}
Let us assume that the dimension $d$ is fixed. The unweighted version of Problem \eqref{mepom} with Euclidean norm can be solved in polynomial time $O(|\A|^{4d+2})$.
\end{thm}
\section{Conclusions\label{sec:7}}
In this paper we analyze minimum radius covering problems with shapes based on polyellipsoids. We study the problem of covering a finite set of $d$-dimensional demand points by optimally finding the dilation and translation of a polyellipsoid with given foci. We provide several formulations for the problem in general dimension $d$: a primal SOC formulation, a conic dual SOC model and a Lagrangean dual program whose resolution needs to solve different Weber problems. We analyze the worst case complexity of these problems and develop an Elzinga-Hearn decomposition approach to solve them geometrically. We report the results of an extensive battery of experiments to show the performance of the proposed approaches. We further analyze the one-dimensional problem and derive a linear time algorithm (in the number of foci) for solving the problem as well as some closed formulas for the covering polyellipse in terms of the foci and the demand points.  Finally we draw some comments on two extensions of the problem: 1) the problem of simultaneously finding the minimum radius and selecting the foci from a potential finite set of candidates; and 2) the use of more general shapes, which we call ordered median polyellipsoids, to cover the demand points. As far as  we know, both extensions, have been first introduced and analyzed here.

Further research on the topic includes, among others, the analysis of the multiple facility case, in which several polyellipsoids have to be located by minimizing the maximum  sum of the closest distances from the demand points to the foci, extending the multifacility center problems; or maximal covering location problems using polyellipsoidal shapes.

\end{document}